\documentclass[11pt]{article}
\usepackage[T1]{fontenc}
\usepackage[utf8]{inputenc}
\usepackage{amsmath,amssymb,amsthm,microtype}
\usepackage{array,longtable}
\usepackage[colorlinks=true,linkcolor=blue,citecolor=blue,urlcolor=blue]{hyperref}
\usepackage[margin=1in]{geometry}
\newtheorem{theorem}{Theorem}[section]
\newtheorem{proposition}[theorem]{Proposition}
\newtheorem{lemma}[theorem]{Lemma}
\newtheorem{corollary}[theorem]{Corollary}
\theoremstyle{definition}
\newtheorem{definition}[theorem]{Definition}
\newtheorem{remark}[theorem]{Remark}

\newtheorem{conjecture}[theorem]{Former conjecture}
\newcommand{\R}{\mathbb R}
\newcommand{\Q}{\mathbb Q}
\newcommand{\Z}{\mathbb Z}
\newcommand{\T}{\mathbb T}
\newcommand{\tu}{\tau_{\mathrm{usual}}}

\newcommand{\Hom}{\operatorname{Hom}}
\newcommand{\init}{\operatorname{init}}
\newcommand{\Comp}{\operatorname{comp}}

\title{Compatible Hamel Topologies, Real Valued Isomorphisms and Rational Scalar Continuity}
\author{Irina Yaroshevskaya}
\date{}
\hypersetup{pdftitle={Compatible Hamel Topologies, Real Valued Isomorphisms and Rational Scalar Continuity},pdfauthor={Irina Yaroshevskaya}}
\begin{document}
\maketitle
\begin{abstract}
This article collects two linked studies of Hamel character topologies on the real line, written as companions of the article on almost periodic characters and Hamel coefficients in TVG topologies. Part I, the more general of the two, treats both the full real-orbit (TVG) dual pair and the rational-parameter (TVGQ) dual pair. It proves a divisibility-chain criterion for compatibility on the rationals and on the torsion group $\Q/\Z$. Compatible subsets of one chain have a compatible join, but two different chains need not: an explicit Chinese remainder construction gives two compatible topologies with an incompatible join, answering the compact-union question negatively and disproving the Greatest Compatible TVG Conjecture for every nonzero rational-valued Hamel functional. The corresponding rational and real Hamel dual pairs admit no Mackey topology. Part I also classifies additive homomorphisms and isomorphisms of the real-valued observation spaces and computes their rational-parameter automorphism groups, showing that nonscalar shears coexist with a scalar criterion for isomorphism types, shows by the Glicksberg theorem for nuclear groups that both prime square-denominator topologies are non-nuclear, classifies sequences converging for every coefficient functional, constructs an infinite-dimensional rational-span sequence whose convergence depends on that functional, and shows that the set of real parameters of continuous characters of the integers with the mixed-prime metric has Hausdorff dimension zero.

Part II studies real vector spaces with Hausdorff additive group topologies admitting every fixed rational scalar multiplication continuously. A rational-valued Hamel functional and a rational parameter space $S$ generate precompact, modular, and real-valued observation topologies, whose exact scalar ring is the multiplier ring of $S$. We compute their completions and modular character groups and classify the precompact TVGQ topologies. A universal modular versus real-valued separating sequence exists exactly when finitely many real scalar multiples of $S$ span the real line; this criterion resolves the case of proper parameter spaces of full rational dimension, where both outcomes occur. Local compactness and sigma-compactness force the continuous scalar ring to be a field; orbit continuity and rational equicontinuity imply joint real scalar continuity. We classify additive isomorphisms within each modular family and obtain exactly $2^{\mathfrak c}$ additive topological group isomorphism types in each rational-parameter modular family. Compatible rational scalar saturation identifies greatest-compatible TVGQ existence with the unrestricted Mackey problem for the same dual pair, and all compatible topologies for the rational Hamel pair are non-Baire. Part II applies the two-chain counterexample of Part I to the rational-scalar framework, where its rational and real-line transfers disprove greatest-compatible existence for every nonzero Hamel functional and rule out Mackey topologies for the associated dual pairs. The operator framework and compact duality are established background; no historical priority claim is made.
\end{abstract}

\tableofcontents

\section*{About this article}
\addcontentsline{toc}{section}{About this article}
This article joins two studies of Hamel character topologies on the real line that accompany \cite{HamelCompanion}. Part I, previously circulated as \emph{Compatible Hamel Topologies and Real Valued Isomorphisms}, is the more general of the two. Its chain criteria, its two-chain counterexample and the transfers of that counterexample cover the full real-orbit (TVG) dual pair $H_f$ as well as the rational-parameter (TVGQ) pair $H_{f,\Q}$. Its isomorphism classification allows arbitrary rational parameter spaces, and its sequence results concern the full real parameter space. Part II, previously \emph{Rational Scalar Continuity and Hamel Character Topologies on the Real Line}, develops the framework for an arbitrary nonzero rational parameter space and applies the two-chain counterexample of Part I to it.

Each part keeps its own section and theorem numbering, exactly as in the separate manuscript; a result of the other part is cited with the name of the part (for example, Theorem 3.3 of Part I). In Part II the references to the former separate preprint are replaced by references to the corresponding results of Part I, and the two bibliographies are merged into one. No statement or proof has been added to either part: apart from these cross-references, the wording that described the parts as separate articles, and the single acknowledgement and bibliography at the end, the text of each part is that of the corresponding separate manuscript. Part I does not depend on Part II, apart from the terminology TVGQ.

\part{Compatible Hamel Topologies and Real Valued Isomorphisms}\label{part:I}

\section{Relation to the preceding article}
The article \cite{HamelCompanion} determines the modular Hamel character groups and their scalar rigidity. Its Conjecture 23.2 asks whether the compatible locally quasi-convex TVG topologies have a greatest member; Question B.4 asks whether compatibility is preserved by unions of profinitely null compact observation sets. Theorem 24.3 there reduces the first question to a suitable pair of topologies on $\Q/\Z$, but does not construct that pair. Appendix C excludes splittings of the dyadic chain.

Here we give the missing pair and answer both questions negatively. Section~\ref{sec:chains-audited} extends the dyadic analysis to general divisibility chains, including an exact compatibility criterion and a factorial null sequence. Section~\ref{sec:twochains} constructs two distinct compatible chains with an incompatible join and gives the transfers to the real line. Section~\ref{sec:tauiso} treats the real-valued observation spaces: nonscalar finite-rank perturbations occur, while existence of an isomorphism still reduces to a scalar criterion. Section~\ref{sec:nuclearity} settles nuclearity of both prime topologies by an application of the Glicksberg theorem for nuclear groups. Section~\ref{sec:sequences} classifies convergence for every coefficient functional and shows that infinite-dimensional convergence can depend on the functional. Section~\ref{sec:mixed-dimension} gives a quantitative covering bound for the possible character parameters of the mixed-prime join and proves that they form a set of Hausdorff dimension zero. These results supplement the character, completion, compactness and modular rigidity theorems of \cite{HamelCompanion}.

This part gives the full proofs of the chain criteria, the two-chain counterexample and its rational and real transfers, and the real-valued homomorphism and isomorphism classification in Sections~\ref{sec:chains-audited}--\ref{sec:tauiso}. These sections serve as the primary source for the applications to rational-scalar Hamel topologies in Part II. The nonexistence conclusions for real-scalar and rational-scalar topologies use the same quotient counterexample.

Write $\T=\R/\Z$, $p(t)=t+\Z$, and $\|p(t)\|_{\T}=\inf_{z\in\Z}|t-z|$. A character is a homomorphism into $\T$. A Hausdorff group topology is locally quasi-convex if it has a zero-neighborhood base of sets separated from each exterior point by continuous characters and the closed quarter arc. Topologies are compatible with a prescribed dual when their full continuous character group equals that dual. A Mackey topology means a greatest compatible Hausdorff locally quasi-convex group topology.

A TVG topology on $\R$ makes every fixed real multiplication continuous; TVGQ requires this only for rational multiplications. Joint continuity of the scalar action is not required. Fix a nonzero rational-linear $f:\R\to\Q$, put $f_a(x)=f(ax)$, and set
\[
\Gamma=\{q\mapsto p(rq):r\in\Q\},\qquad
H_f=\{x\mapsto p(cx+f(ax)):c,a\in\R\},
\]
\[
H_{f,\Q}=\{x\mapsto p(cx+qf(x)):c\in\R,\ q\in\Q\}.
\]
We use the quarter-arc criterion: a character is continuous if it maps some zero-neighborhood into $[-1/4,1/4]+\Z$. For completeness, if that bound holds on $U$, choose a neighborhood whose first $N$ integer multiples lie in $U$; Lemma~\ref{chains:small} then bounds the character there by $1/(4N)$.

\section{Divisibility chains and compatible joins}\label{sec:chains-audited}
The following arguments extend the dyadic obstruction. All chains in this section are chains of positive integers, with each term dividing the next and each consecutive ratio at least two. The criteria concern subsets of one chain. They do not decide compatibility of arbitrary joins, and no historical priority is asserted.

Write $p(t)=t+\Z$ and $\|p(t)\|_{\T}=\inf_{z\in\Z}|t-z|$. For an increasing divisibility chain $(M_i)$ put $\alpha_i=M_{i+1}/M_i$. We say that $H(M)$ holds if, for every $A>2$, there is an integer $L(A)\ge0$ such that every run of consecutive steps satisfying $\alpha_i<A$ has at most $L(A)$ steps. For a finite chain this condition holds automatically. If $U$ is a subset of the indices of a fixed chain $(m_n)$, $H(U)$ means this condition for the increasing enumeration $(m_{u_i})$ of the selected terms. Thus the condition concerns ratios of consecutive selected terms, not just consecutive original indices.

\begin{lemma}[Small multiples in a circle arc]\label{chains:small}
If $t\in\T$ and $\|st\|_{\T}\le\frac14$ for $1\le s\le N$, then $\|t\|_{\T}\le\frac{1}{4N}$.
\end{lemma}
\begin{proof}
Choose the centered representative $y$ of $t$. The case $s=1$ gives $|y|\le\frac14$. If $N|y|>\frac14$, let $s\le N$ be the first positive integer with $s|y|>\frac14$. Then $\frac14<s|y|\le\frac12$, contradicting the assumed circle bound. Thus $N|y|\le\frac14$.
\end{proof}

\begin{lemma}[Finite unions of chain subsets]\label{chains:union}
If nonempty subsets $S,T$ of the indices of one divisibility chain satisfy $H(S)$ and $H(T)$, then $H(S\cup T)$ holds.
\end{lemma}
\begin{proof}
Suppose instead that the selected union has arbitrarily long runs with consecutive ratios less than some $A>2$. Put $l=L_T(A)$, $g=l+2$, and $B=A^g$. In such a run a block containing no $S$-term consists of consecutive $T$-terms, so has at most $l+1$ terms. This applies to both end blocks as well as the blocks between successive $S$-terms. If the run has $r$ terms belonging to $S$, its total number of terms is at most
\[
r+(r+1)(l+1)=(r+1)(l+2)-1.
\]
Consequently $r$ can be arbitrarily large. Successive $S$-terms in this run are consecutive in the global selected $S$-chain and have at most $g$ union steps between them. Their ratios are less than $B$. Taking $r\ge L_S(B)+2$ contradicts $H(S)$. The proof also covers finite $S$ or $T$.
\end{proof}

\textbf{Rational observation topologies.}
Fix a chain $(m_n)$ and a nonempty index set $U$. On the additive group $\Q$ put
\[
Q_{U,d}(q)=\sup_{n\in U}\left\|p\left(\frac{q}{d m_n}\right)\right\|_{\T},\qquad d\ge1,
\]
and let $\xi_U$ have zero-neighborhood base $\{q:Q_{U,d}(q)<\varepsilon\}$, with $d$ a positive integer and $\varepsilon>0$. Put $\Gamma=\{q\mapsto p(rq):r\in\Q\}$. These are Hausdorff metrizable locally quasi-convex group topologies admitting every fixed rational scalar continuously. Indeed circle distance is subadditive, and a common multiple $D$ of finitely many denominator parameters gives $Q_{U,d}\le(D/d)Q_{U,D}$. Moreover
\[
Q_{U,d}\left(\frac{a}{b}q\right)\le |a|Q_{U,bd}(q),\qquad a\in\Z,\quad b\ge1.
\]
Their defining characters generate $\Gamma$, which separates points. Intersections of character preimages of closed arcs of radius $1/(4N)$ give quasi-convex neighborhoods: these arcs are the intersection of the inverse images of the closed quarter arc under multiplication by $1,\ldots,N$, by Lemma~\ref{chains:small}. This proves local quasi-convexity. The denominator parameters and a countable choice of radii give metrizability. Each $\xi_U$ is coarser than usual. Therefore every continuous character is $q\mapsto p(cq)$ for some real $c$, by extension of usual-continuous characters from $\Q$ to its usual completion $\R$.

\begin{theorem}[Rational chain criterion]\label{chains:F}
For every chain $(m_n)$ and nonempty $U$, the full character group of $\xi_U$ is $\Gamma$ if and only if $H(U)$ holds. If $H(U)$ fails, an irrational-parameter character is continuous.
\end{theorem}
\begin{proof}
Write the selected terms as $M_i=m_{u_i}$, with ratios $\alpha_i$. First suppose $H(U)$ fails. Choose $A>2$ admitting arbitrarily long runs of ratios less than $A$, and put $\delta=1/(2(A+1))$. Choose disjoint increasing runs $[a_s,b_s]$, $s\ge1$, of at least $s+1$ steps, so far apart that
\[
\frac{M_{b_{s+1}}}{M_{b_s}}\ge2^{s+1}.
\]
Such runs exist arbitrarily far out: deleting finitely many indices cannot bound the lengths of all remaining runs if their lengths were unbounded originally. Put $\gamma=\sum_{s\ge1}1/M_{b_s}$. If $Q_{U,1}(q)<\delta$, let $y_i$ be the centered representative of $p(q/M_i)$. On each selected run, $y_i-\alpha_i y_{i+1}$ is an integer of absolute value less than $(A+1)\delta<1$, hence zero. Thus
\[
|y_{b_s}|<\delta\,2^{-(s+1)}.
\]
For fixed $q$, one has $y_{b_s}=q/M_{b_s}$ eventually. Hence $p(\gamma q)=p(\sum_s y_{b_s})$, and $\|p(\gamma q)\|_{\T}<\delta/2<1/4$. The quarter-arc criterion gives continuity. For irrationality, the $k$th partial sum has denominator dividing $M_{b_k}$, while its strictly positive tail is at most $2/M_{b_{k+1}}$. If $\gamma=u/v$ were rational, multiplying that tail by $vM_{b_k}$ would give a positive integer smaller than one for large $k$. This is impossible.

Now suppose $H(U)$ holds. If $U$ is finite, $\xi_U$ is the weak topology of its defining rational characters, whose generated character group is exactly $\Gamma$; precompact duality gives the assertion. Assume henceforth that $U$ is infinite and $q\mapsto p(cq)$ is continuous. Choose an integer $d\ge1$ and $0<\varepsilon<1/8$ so that
\[
Q_{U,d}(q)<\varepsilon\quad\Longrightarrow\quad\|p(cq)\|_{\T}<1/4.
\]
Put $c'=dc$. Substituting $q=dq_0$ gives this implication with $Q_{U,1}(q_0)$ and $p(c'q_0)$, in particular for integer $q_0$. Fix
\[
h=\frac{\varepsilon}{2},\qquad \kappa=\frac{1}{\varepsilon},\qquad
\eta=\frac{\varepsilon}{16},\qquad B=\frac{64}{\varepsilon},
\]
and put $G=\{i:\alpha_i\ge B\}$. This set is infinite, and successive members $i<k$ have at most $L=L_U(B)$ intervening steps with ratios less than $B$.

Let $r_i$ be the centered representative of $p(c'M_i)$. For $i\in G$ and $1\le s\le\lfloor h\alpha_i\rfloor$, the integer $sM_i$ satisfies
\[
Q_{U,1}(sM_i)\le s/\alpha_i\le h<\varepsilon.
\]
Lemma~\ref{chains:small}, and $\lfloor h\alpha_i\rfloor\ge h\alpha_i/2$, give
\[
|r_i|\le\frac{\kappa}{\alpha_i}\le\frac{1}{64}.
\]
For successive $i<k$ in $G$, write $D_i=M_k/M_{i+1}\le B^L$ and $P_i=M_k/M_i=\alpha_iD_i$. Then
\[
b_i=P_i r_i-r_k\in\Z.
\]
If $b_i=0$ eventually along $G$, the relation $r_k=P_i r_i$ multiplies a nonzero absolute value by at least two at every step, contradicting the uniform bound above. Thus these $r_i$ are zero on a tail of $G$, and $c'M_i\in\Z$ for some $i$. Hence $c\in\Q$.

Otherwise $b_i\ne0$ for infinitely many such pairs. At their first indices,
\[
 c_0\le\alpha_i|r_i|\le\kappa,\qquad c_0=\frac{3}{4B^L}>0,
\]
because $|P_i r_i|\ge1-|r_k|\ge3/4$. Enumerate these indices increasingly as $i_j$, put $s_j=\lfloor\eta\alpha_{i_j}\rfloor$ and $\sigma_j=\operatorname{sgn}(r_{i_j})$, and consider the integer
\[
q_0=\sum_{j=1}^{m}\sigma_j s_jM_{i_j}.
\]
For any selected index $l$, terms with $i_j\ge l$ vanish modulo $M_l$. The others have absolute contributions at most $\eta M_{i_j+1}/M_l$. Distinct such indices give a geometric sum at most two. Therefore
\[
Q_{U,1}(q_0)\le2\eta<\varepsilon.
\]
On the other hand $p(c'q_0)=p(\sum_{j=1}^{m}s_j|r_{i_j}|)$, and every summand belongs to $[\eta c_0/2,\eta\kappa]$, where $\eta\kappa=1/16$. Choose the first $m$ for which this positive partial sum exceeds $1/4$. Its value is at most $5/16$, so its circle norm exceeds $1/4$, a contradiction. Only the rational case is possible, proving the asserted full dual.
\end{proof}

\begin{corollary}[No incompatible pair from subsets of one chain]\label{chains:Fjoin}
For nonempty index sets $S,T$ of one chain,
\[
\xi_S\vee\xi_T=\xi_{S\cup T}.
\]
If both topologies have full dual $\Gamma$, so does their join.
\end{corollary}
\begin{proof}
For the same denominator parameter, $Q_{S\cup T,d}=\max\{Q_{S,d},Q_{T,d}\}$. Common multiples of denominator parameters handle every finite intersection, giving the equality. Apply Theorem~\ref{chains:F} and Lemma~\ref{chains:union}.
\end{proof}

\begin{proposition}[Factorial chains and an explicit null sequence]\label{chains:factorial}
For the chain $1!,2!,3!,\ldots$, all nonempty subset topologies above have full dual $\Gamma$. The full-chain topology is strictly coarser than the usual topology. A sequence tending to zero in it is
\[
N_j=(j!)^2,\qquad q_j=j!\,N_j!,
\]
although $q_j$ tends to infinity usually.
\end{proposition}
\begin{proof}
The consecutive ratios tend to infinity, so every selected subsequence satisfies $H$. Fix $d$. Eventually $d\mid j!$. For $k\le N_j$, $q_j/(dk!)$ is an integer. For $k>N_j$ its absolute value is at most
\[
\frac{j!}{d(N_j+1)}<\frac{1}{d\,j!}\longrightarrow0.
\]
This proves convergence for every denominator parameter. The denominator parameter $d$ is essential: $(n!)!$ does not converge to zero, because the observation $k=n!$ with $d=2$ has circle norm $1/2$.
\end{proof}

\textbf{Observation topologies on the torsion quotient.}
Put $A=\Q/\Z$ and let $\iota:A\to\T$ be the natural inclusion. Every algebraic character on $A$ is identified with a profinite integer $\zeta\in\widehat{\Z}$: its values at $1/b+\Z$ give consistent residues modulo every positive integer $b$. For $C\subseteq\Z$ containing $1$ put
\[
P_C(x)=\sup_{c\in C}\|cx\|_{\T},\qquad
\lambda_C:\quad\{x:P_C(x)<\varepsilon\},\quad\varepsilon>0.
\]
These are Hausdorff locally quasi-convex group topologies by the same arc argument as above. Compatibility means that the full dual is exactly $\Z\iota$. Adding $0$ to $C$ does not change $P_C$.

\begin{theorem}[Torsion chain criterion]\label{chains:G}
Let $C$ be a finite divisibility chain containing $1$, or an infinite chain $C=\{M_0,M_1,\ldots\}$ with $M_0=1$ and ratios $\alpha_i\ge2$. In the infinite case assume that every positive integer divides some $M_i$. Then $\lambda_C$ is compatible if and only if $H(C)$ holds.
\end{theorem}
\begin{proof}
For finite $C$, $\lambda_C$ is the circle subspace topology: integer multiplications are circle-continuous and the observation $1$ is included. Its dual is $\Z\iota$. Assume $C$ infinite.

If $H(C)$ fails, choose $A_0>2$ with arbitrarily long runs of ratios less than $A_0$, and $0<\delta<1/(2(A_0+1))$. Choose runs $[a_s,b_s]$ of at least $s+1$ steps, with $a_{s+1}\ge b_s+2$. For $P_C(x)<\delta$ use the centered representatives $y_i$ of $M_ix$, not their fractional parts. The integer $y_{i+1}-\alpha_i y_i$ has absolute value less than $(A_0+1)\delta<1$ on a run, hence is zero. Consequently
\[
|y_{a_s}|<\delta\,2^{-(s+1)}.
\]
The series $\zeta=\sum_{s\ge1}M_{a_s}$ converges in $\widehat{\Z}$, because every modulus divides all sufficiently late terms. It is not an ordinary integer: its mixed-radix digits are $1$ at the indices $a_s$ and $0$ elsewhere, with infinitely many nonzero digits and infinitely many intervening zero digits. A nonnegative integer has eventually zero digits; a negative integer has eventually the maximal digits $\alpha_i-1$. Neither is possible here. For a fixed $x\in A$ all sufficiently late terms $M_{a_s}x$ vanish. Thus
\[
\|\zeta x\|_{\T}\le\sum_{s\ge1}|y_{a_s}|<\delta/2<1/4.
\]
The quarter-arc criterion makes this noninteger character continuous, proving incompatibility.

Suppose now that $H(C)$ holds and an algebraic character $\zeta$ is continuous. Choose $0<\varepsilon<1/8$ so that $P_C(x)<\varepsilon$ implies $\|\zeta x\|_{\T}<1/4$. Use the same constants
\[
h=\frac{\varepsilon}{2},\qquad \kappa=\frac{1}{\varepsilon},\qquad
\eta=\frac{\varepsilon}{16},\qquad B=\frac{64}{\varepsilon},
\]
and $G=\{i:\alpha_i\ge B\}$, with at most $L=L_C(B)$ consecutive steps outside $G$. Write $\theta_i=z_i/M_{i+1}$, where $z_i$ is the residue of $\zeta$ in $\{0,\ldots,M_{i+1}-1\}$, and let $r_i$ be the centered representative of $p(\theta_i)$.

For $i\in G$ and $1\le s\le\lfloor h\alpha_i\rfloor$, the point $x=s/M_{i+1}+\Z$ has $P_C(x)\le h<\varepsilon$: observations above $i$ are integers, and observations at or below $i$ have absolute representatives at most $s/\alpha_i$. Lemma~\ref{chains:small} therefore gives
\[
|r_i|\le\frac{\kappa}{\alpha_i}\le\frac{1}{64}.
\]
For successive $i<k$ in $G$, put $D_i=M_k/M_{i+1}\le B^L$ and $P_i=M_{k+1}/M_{i+1}=\alpha_kD_i$. Consistency of residues gives
\[
b_i=P_i r_k-r_i\in\Z.
\]
If $b_i=0$ eventually, the signed integers $u_i=M_{i+1}r_i$ are constant along a tail of $G$. Let the common value be $u$. Then $\zeta\equiv u$ modulo each of the cofinal moduli $M_{i+1}$ on that tail, so $\zeta=u$ in $\widehat{\Z}$. This includes negative integers: their positive residues grow with the modulus, while these signed representatives stay constant.

Otherwise $b_i\ne0$ for infinitely many pairs. At their second indices $k$,
\[
c_0\le\alpha_k|r_k|\le\kappa,\qquad c_0=\frac{3}{4B^L}>0.
\]
Enumerate those indices increasingly as $k_j$, put $s_j=\lfloor\eta\alpha_{k_j}\rfloor$ and $\sigma_j=\operatorname{sgn}(r_{k_j})$, and consider
\[
x=\sum_{j=1}^{m}\frac{\sigma_js_j}{M_{k_j+1}}+\Z.
\]
For any $l$, terms with $k_j<l$ become integers after multiplication by $M_l$. The others have absolute contributions at most $\eta M_l/M_{k_j}$, whose sum is at most $2\eta$. Therefore $P_C(x)\le2\eta<\varepsilon$. But
\[
\zeta x=p\left(\sum_{j=1}^{m}s_j|r_{k_j}|\right),\qquad
\frac{\eta c_0}{2}\le s_j|r_{k_j}|\le\eta\kappa=\frac{1}{16}.
\]
The first positive partial sum exceeding $1/4$ is at most $5/16$, contradicting the continuity bound. Thus $\zeta$ must be an integer, proving compatibility.
\end{proof}

\begin{corollary}[The compact-union question inside one chain]\label{chains:Gjoin}
Let $C_1,C_2$ be subsets of one infinite divisibility chain starting with $1$, such that every positive integer divides some term. Suppose both subsets contain $1$. If $\lambda_{C_1}$ and $\lambda_{C_2}$ are compatible, then
\[
\lambda_{C_1\cup C_2}=\lambda_{C_1}\vee\lambda_{C_2}
\]
is compatible. Finite subsets are allowed.
\end{corollary}
\begin{proof}
Each infinite subset inherits the divisibility-cofinality assumption. Theorem~\ref{chains:G} and Lemma~\ref{chains:union} prove the result for two infinite subsets. A finite subset gives the circle topology, which is coarser than every $\lambda_C$ with $1\in C$, so the finite cases follow as well. The equality of joins follows from $P_{C_1\cup C_2}=\max\{P_{C_1},P_{C_2}\}$. Each infinite subset, together with $0$, is compact in the profinite topology on $\Z$, because its terms converge to zero profinitely. Thus this is indeed the stated restricted compact-union problem.
\end{proof}

The same-chain restriction is essential: Theorem~\ref{twochains:H} gives different chains answering Question B.4 of \cite{HamelCompanion} negatively. The rational and torsion criteria concern their respective constructions; no identification of their topology lattices is used.

\textbf{Irrational scalar probes for the rational Hamel pair.}
Here a TVGQ topology means a Hausdorff additive group topology admitting every fixed rational multiplication continuously. This paragraph concerns the rational-parameter pair, not the full real-orbit character group. Fix a nonzero rational-valued additive $f:\R\to\Q$ and write
\[
H_{f,\Q}=\{x\mapsto p(cx+qf(x)):c\in\R,\ q\in\Q\}.
\]
For comparison, let $\beta_{f,\Q}$ be the initial topology of the characters in $H_{f,\Q}$, put $\sigma_{f,\Q}=\tu\vee\beta_{f,\Q}$, and let $\tau_{f,\Q}$ be the initial topology of the usual identity coordinate and $f:\R\to\Q_{\mathrm{usual}}$. These definitions restrict the Hamel observations to rational parameters.

\begin{proposition}[Excluded scalars and sequence witnesses]\label{chains:irrational}
For every irrational $r$, the character $\psi_r(x)=p(f(rx))$ is outside $H_{f,\Q}$. Every compatible TVGQ topology with this prescribed dual has continuous scalar ring exactly $\Q$. For each fixed $f$ and irrational $r$, there is a sequence $(x_n)$ with
\[
|x_n|<1/n,\qquad f(x_n)=0,\qquad f(rx_n)=\frac12.
\]
It converges to zero in $\beta_{f,\Q}$, $\sigma_{f,\Q}$, and $\tau_{f,\Q}$, whereas $(rx_n)$ converges to zero in none of these three topologies.
\end{proposition}
\begin{proof}
If $\psi_r\in H_{f,\Q}$, equality of the circle characters gives $f(rx)-cx-qf(x)\in\Z$ for all $x$. This additive integer-valued function vanishes, by divisibility of the additive real group. The rational-valued difference $f(rx)-qf(x)$ then equals $cx$, forcing $c=0$. Thus $f((r-q)x)=0$ for all $x$, which contradicts $f\ne0$ and $r\notin\Q$. If multiplication by $r$ were continuous in a compatible topology, composition with the prescribed character $p\circ f$ would make $\psi_r$ continuous, a contradiction. Rational multiplications are continuous by the TVGQ assumption, proving the scalar-ring assertion.

The two rational-linear maps $f$ and $x\mapsto f(rx)$ are independent. Indeed a nontrivial rational linear dependence would make $f((a+br)x)=0$ for rational $a,b$ not both zero, impossible for irrational $r$. Their joint map is onto $\Q^2$, and its kernel has finite rational codimension in $\R$. This kernel is nonzero and contains a usual-dense rational line. Thus the fiber over $(0,1/2)$ meets $(-1/n,1/n)$ for each $n$, giving the sequence. Its usual coordinate tends to zero and every rational multiple of $f(x_n)$ vanishes, proving convergence in all three initial topologies. But the character $p\circ f$, continuous in all three, takes the constant value $1/2+\Z$ on $rx_n$. Hence no claimed image convergence holds. The sequence is chosen after $f$ and $r$; no common sequence for all functionals, or convergence in every finer compatible topology, is asserted.
\end{proof}

\begin{proposition}[A fixed irrational character as a join test]\label{chains:probe}
Fix an irrational $r$. Suppose $\eta_1,\eta_2$ are compatible locally quasi-convex TVGQ topologies with full dual $H_{f,\Q}$. If respective zero-neighborhoods $U_1,U_2$ satisfy
\[
\psi_r(U_1\cap U_2)\subseteq[-\tfrac14,\tfrac14]+\Z,
\]
then their join is incompatible. Thus such a pair would rule out a greatest compatible member for this rational Hamel pair.
\end{proposition}
\begin{proof}
The intersection is a zero-neighborhood in the join. The quarter-arc criterion makes $\psi_r$ continuous there, while Proposition~\ref{chains:irrational} places it outside the prescribed dual. A greatest compatible member could not dominate that join. The hypotheses do not construct the two topologies or the required neighborhoods.
\end{proof}

The continuous scalars of any TVGQ form a $\Q$-subalgebra $C_\tau$ of $\R$. If $r\notin C_\tau$, then $ar+b\notin C_\tau$ for rational $a\ne0,b$, since otherwise $r=(ar+b-b)/a$ would belong to the ring. Thus excluded scalars are dense whenever the topology is not a real TVG. The algebraic elements in $C_\tau$ form a field: the inverse of a nonzero algebraic element is a rational polynomial in that element. Consequently algebraic $r\ne0$ and $1/r$ have the same continuity status. This need not hold for transcendental scalars when the continuous scalar ring is not a field.

The probe above is specific to the rational parameter restriction. For the full real-orbit pair $H_f=\{x\mapsto p(cx+f(ax)):c,a\in\R\}$, every $\psi_r$ already belongs to $H_f$ and cannot witness an extra character. The chain results therefore provide restricted exclusions, and the irrational probe provides a conditional counterexample test for TVGQ; the separate two-chain construction in Section~\ref{sec:twochains} supplies the counterexamples to both greatest-compatible assertions.

\section{Two chains with an incompatible join}\label{sec:twochains}
The single-chain union result does not extend to two different chains. We construct two compatible topologies whose join has an extra character. All compatibility assertions use the full character group, and all scalar maps required for the two real-line transfers are included explicitly.

\begin{theorem}[Two compatible chains with no compatible common upper bound]\label{twochains:H}
There are two divisibility chains $F_1,F_2$ of positive integers, each containing $1$, with consecutive ratios tending to infinity and with every positive integer dividing a term of each chain, such that $\lambda_{F_1}$ and $\lambda_{F_2}$ have full dual $\Z\iota$, while their join has a continuous character in $\widehat{\Z}\setminus\Z$. The sets $C_i=\{0\}\cup F_i$ are compact for the profinite topology on $\Z$.
\end{theorem}
\begin{proof}
First construct $g_0=1\mid g_1\mid\cdots$ in stretches indexed by $s=1,2,\ldots$. Stretch $s$ consists of $s+3$ steps of ratio $2$, followed by one step of ratio $s+1$. Let $a_s$ be the starting index of its ratio-$2$ run and $b_s=a_s+s+3$. At the start of stretch $s$, the product of the preceding exceptional ratios is $s!$. Hence $s!\mid g_{a_s}$, and every positive integer divides all sufficiently late $g_n$. The chain has arbitrarily long ratio-$2$ runs. Put
\[
t_n=\max\{t\ge1:t!\mid g_n\}.
\]
These finite integers are nondecreasing and tend to infinity.

Set $M_0=1$ and $E_n=M_n+g_n$, so $E_0=2$. Given $M_n$, put $L_n=\operatorname{lcm}(M_n,t_{n+1}!)$ and require
\[
M_{n+1}\equiv0\pmod{L_n},\qquad
M_{n+1}\equiv-g_{n+1}\pmod{E_n}.
\]
The Chinese remainder theorem applies. Indeed,
\[
\gcd(M_n,E_n)=\gcd(M_n,g_n)\mid g_{n+1},\qquad
\gcd(t_{n+1}!,E_n)\mid g_{n+1}.
\]
For positive integers $a,b,c$, prime valuations give
\[
\gcd(\operatorname{lcm}(a,b),c)
=\operatorname{lcm}(\gcd(a,c),\gcd(b,c)).
\]
Therefore $\gcd(L_n,E_n)\mid g_{n+1}$, exactly the compatibility condition for the two congruences. Their solutions form an arithmetic progression with positive modulus, so choose a positive solution satisfying
\[
M_{n+1}\ge\max\{(n+2)^2M_n,\;2^{n+1}g_{n+1}\}.
\]
This proves $M_n\mid M_{n+1}$ and $E_n\mid E_{n+1}$. Inductively $M_n\ge2^n g_n$, so $E_n\le2M_n$ and
\[
\frac{M_{n+1}}{M_n}\ge(n+2)^2,\qquad
\frac{E_{n+1}}{E_n}\ge\frac{(n+2)^2}{2}.
\]
Both ratios are at least two and tend to infinity. Also $t_n!\mid M_n$ and $t_n!\mid E_n$ for every $n$. Consequently
\[
F_1=\{M_n:n\ge0\},\qquad F_2=\{1\}\cup\{E_n:n\ge0\}
\]
meet all the hypotheses of Theorem~\ref{chains:G}. Ratios tending to infinity imply $H(F_i)$: for each fixed bound only finitely many ratios fall below it. Thus each $\lambda_{F_i}$ has full dual $\Z\iota$.

Write $G=\{g_n:n\ge0\}$. The difference identity $E_n-M_n=g_n$ gives, for every $x\in\Q/\Z$,
\[
P_G(x)\le P_{F_1}(x)+P_{F_2}(x).
\]
It follows that the join $\lambda_{F_1}\vee\lambda_{F_2}$ is finer than $\lambda_G$. The failure direction of Theorem~\ref{chains:G} already supplies an extra character. More explicitly, let
\[
\zeta_G=\sum_{s\ge1}g_{a_s}\in\widehat{\Z}.
\]
The series converges because its terms tend to zero modulo every positive integer. Its mixed-radix digits in the $g_n$ chain are $1$ at the indices $a_s$ and $0$ elsewhere. There are infinitely many nonzero digits and infinitely many zero digits. It is neither a nonnegative integer, whose digits are eventually zero, nor a negative integer, whose digits are eventually maximal. Thus $\zeta_G\notin\Z$.

For a direct continuity check fix $0<\delta<1/16$ and suppose $P_{F_1}(x)<\delta$ and $P_{F_2}(x)<\delta$. All centered representatives $y_n$ of $g_nx$ have absolute value less than $2\delta$. Along a ratio-$2$ run, $y_{n+1}-2y_n$ is an integer of absolute value less than $6\delta<1$, so it is zero. Hence
\[
|y_{a_s}|<2\delta\,2^{-(s+3)},\qquad
\|\zeta_Gx\|_{\T}\le\sum_{s\ge1}|y_{a_s}|<\frac{\delta}{4}<\frac14.
\]
For each fixed torsion element $x$, all sufficiently late $g_{a_s}x$ vanish, so this sum represents the character value. The quarter-arc criterion proves continuity in the join. A compatible common upper bound would make this same extra character continuous, which is impossible. Finally each of the two chains tends to zero profinitely; adjoining its limit gives a compact set $C_i$, and adjoining $0$ does not change its topology.
\end{proof}

\begin{corollary}[The general compact-union assertion is false]\label{twochains:compact}
Two profinitely compact subsets $C_1,C_2\subseteq\Z$, each containing $0$ and $1$, can define compatible topologies on $\Q/\Z$ while $\lambda_{C_1\cup C_2}$ is incompatible. The positive union result of Corollary~\ref{chains:Gjoin} remains valid when both sets are subsets of one divisibility chain.
\end{corollary}
\begin{proof}
Use the sets in Theorem~\ref{twochains:H} and the identity $P_{C_1\cup C_2}=\max\{P_{C_1},P_{C_2}\}$.
\end{proof}

\begin{theorem}[Consequences for rational and real scalar topologies]\label{twochains:transfer}
For every nonzero rational-valued additive $f:\R\to\Q$, there are two compatible Hausdorff locally quasi-convex TVG topologies for $H_f$ with incompatible join, and two compatible Hausdorff locally quasi-convex TVGQ topologies for $H_{f,\Q}$ with incompatible join. Thus both greatest-compatible assertions are false for every such $f$. None of the dual pairs
\[
(\Q/\Z,\Z\iota),\qquad(\Q,\Gamma),\qquad
(\R,H_f),\qquad(\R,H_{f,\Q})
\]
admits a Mackey topology in the class of Hausdorff locally quasi-convex group topologies.
\end{theorem}
\begin{proof}
We give the transfer details to distinguish the two real-line constructions. Put $\pi(q)=q+\Z$ and let $\lambda_i=\lambda_{F_i}$. Define $\delta_i$ on $\Q$ initially by all maps $q\mapsto\pi(rq)$ into $(\Q/\Z,\lambda_i)$, with $r\in\Q$. Their kernels separate points, the target topologies are locally quasi-convex, and rational rescalings permute the defining maps. Thus $\delta_i$ is Hausdorff, locally quasi-convex, and admits every fixed rational scalar continuously.

Its full dual is $\Gamma$. To verify this, a continuous character has a quarter-arc bound on a neighborhood using finitely many parameters $k_j/d$. Integer multiplication on $(\Q/\Z,\lambda_i)$ is continuous, so those observations factor continuously through the single surjection $h_d(q)=\pi(q/d)$. The same bound makes the character continuous for the pullback topology of $h_d$. It kills $d\Z$, the closure of zero in that pullback topology, and descends to a $\lambda_i$-continuous character. That character is an integer multiple of $\iota$, giving $q\mapsto p(nq/d)$. Conversely every rational-parameter character is continuous. This proves the full dual assertion.

The character $\psi=\zeta_G\circ\pi$ is continuous on $\delta_1\vee\delta_2$: the map $\pi$ is continuous between the corresponding joins. It is outside $\Gamma$. Indeed if $\psi(q)=p(rq)$ with $r\in\Q$, evaluation at $1$ gives $r\in\Z$; surjectivity of $\pi$ would then give $\zeta_G=r\iota$, contradicting Theorem~\ref{twochains:H}.

For the rational-scalar construction put
\[
S_f(\delta)=\operatorname{init}\{\mathrm{id}:\R\to\R_{\mathrm{usual}},\quad f:\R\to(\Q,\delta)\}.
\]
Its scalar maps for rational scalars are continuous. The group $N=\ker f$ is usual-dense, and every fiber of $f$ meets every usual interval. Thus $f:S_f(\delta)\to(\Q,\delta)$ is an open continuous surjection, and the induced topology on $N$ is usual. Restrict a continuous character to $N$ and extend it to the usual completion $\R$, obtaining $p(cx)$. After subtracting that ordinary character it descends through $f$ to a continuous character of $(\Q,\delta)$. If the latter dual is $\Gamma$, this proves that the full dual of $S_f(\delta)$ is exactly $H_{f,\Q}$.

For real scalars put $f_a(x)=f(ax)$ and
\[
R_f(\delta)=\operatorname{init}\{\mathrm{id}:\R\to\R_{\mathrm{usual}},\quad f_a:\R\to(\Q,\delta)\ (a\in\R)\}.
\]
Fixed real dilations permute the Hamel coordinates, giving the TVG property. To check its dual, a character is controlled by finitely many coordinates and the usual coordinate. Choose a rational basis $h_1,\ldots,h_m$ of those Hamel coordinates. Rational coordinate changes are continuous on $(\Q,\delta)^m$. The joint map $h:\R\to\Q^m$ is onto by independence of the functionals, and its kernel is usual-dense because it has finite rational codimension in $\R$. Every fiber meets every usual interval, so $h$ is open for the topology generated by these finitely many coordinates and the usual coordinate. The character is continuous for that finite-coordinate topology by the quarter-arc criterion. Restriction to its kernel and extension to $\R$ again gives an ordinary character $p(cx)$. Subtracting it yields a continuous character of $(\Q,\delta)^m$, a sum of rational-parameter coordinate characters. Every $h_j$ is a rational combination of maps $f_a$, so the result is $p(cx+f(bx))$. This proves that the full dual is exactly $H_f$.

Both initial constructions preserve binary joins and are Hausdorff and locally quasi-convex. In either transferred join the character $\psi\circ f$ is continuous. It lies outside both prescribed duals. For the rational pair, an equality $\psi(f(x))=p(cx+rf(x))$ forces $c=0$ by restriction to the nonzero divisible group $N$, then forces $\psi(q)=p(rq)$ by surjectivity. For the full real-orbit pair, an equality $\psi(f(x))=p(cx+f(ax))$ makes $cx+f(ax)$ integer-valued on $N$. Divisibility makes this additive integer-valued map zero. Hence $cN\subseteq\Q$; since $N$ is uncountable, $c=0$. Now $f_a$ vanishes on $N$, so $f_a=rf$ for some rational $r$, again forcing $\psi\in\Gamma$. Both conclusions contradict the choice of $\psi$.

One can also exclude the extra character from both prescribed duals directly. Choose $u\in\R$ with $f(u)=1$. If $\zeta_G\circ\pi\circ f=p\circ h$ for any real-valued additive $h$, restriction to $u\Q$ gives $(\zeta_G\circ\pi)(q)=p(qh(u))$. Evaluation at $q=1$ forces $h(u)\in\Z$, and surjectivity of $\pi$ then makes $\zeta_G$ an ordinary integer, a contradiction. Every prescribed Hamel character has such a real-valued additive lift.

Thus each pair has compatible members and an incompatible join, even with the required scalar continuity. No compatible topology, including any unrestricted Mackey topology, can dominate either pair. The same argument applies directly to the pairs on $\Q/\Z$ and $\Q$. This proves all the asserted nonexistence statements.
\end{proof}

This construction is independent of the choice of nonzero $f$. It does not determine the character group of the specific mixed-prime join. Nuclearity of its two members is settled in Section~\ref{sec:nuclearity}. Section~\ref{sec:tauiso} separately classifies additive isomorphisms of the real-valued observation spaces. 

\section{Additive isomorphisms of the real-valued observation spaces}\label{sec:tauiso}
Let $f,g:\R\to\Q$ be nonzero rational-linear functionals and $S,T\subseteq\R$ rational subspaces. Write $f_a(x)=f(ax)$ and
\[
\mathcal F(f,S)=\{f_a:a\in S\},\qquad
\tau_{f,S}=\operatorname{init}\{\mathrm{id}:\R\to\R_{\mathrm{usual}},\quad f_a:\R\to\R_{\mathrm{usual}}\ (a\in S)\}.
\]
For $S=\R$ this is the full real-valued topology $\tau_f$. The modular rigidity result does not extend to all homomorphisms between these spaces. Nevertheless their isomorphism types have the same scalar criterion. Empty observation spaces are allowed in the structural statements below.

\begin{lemma}[Continuous real-valued and rational-valued maps]\label{tauiso:maps}
If $B$ is a rational basis of $S$, the continuous real-valued additive maps on $(\R,\tau_{f,S})$ are exactly
\[
h(x)=cx+\sum_{j=1}^n d_jf(b_jx),\qquad c,d_j\in\R,\quad b_j\in B.
\]
Its continuous rational-valued additive maps are exactly $\mathcal F(f,S)$.
\end{lemma}
\begin{proof}
Distinct basis elements give independent rational-valued functionals: a rational relation is $f_a=0$, which forces $a=0$, since multiplication by nonzero $a$ is onto. The joint map of any finite independent list is onto $\Q^n$. Its kernel has finite rational codimension, is uncountable, and is usual-dense because it contains a nonzero rational line. Every fiber therefore meets every usual interval. It follows that the graph of the usual coordinate and this finite list is dense in $\R^{n+1}$.

A continuous additive map into the usual real line is controlled by finitely many defining coordinates: a bound on one basic neighborhood gives continuity for that finite-coordinate topology by taking smaller neighborhoods whose first finitely many integer multiples lie in the bounded one. It extends to its completion $\R^{n+1}$. A continuous additive map on this finite-dimensional real space is real-linear, giving the displayed formula. Conversely each such formula is continuous by definition.

If $h$ is rational-valued, its restriction to the common kernel gives $cx\in\Q$ on an uncountable group, so $c=0$. Surjectivity of the joint map onto $\Q^n$, tested at the rational standard unit vectors, gives $d_j\in\Q$. Thus $h=f_{\sum_jd_jb_j}\in\mathcal F(f,S)$. The empty-list case gives only the zero rational-valued map.
\end{proof}

\begin{theorem}[All continuous additive homomorphisms]\label{tauiso:hom}
An additive map $L:(\R,\tau_{f,S})\to(\R,\tau_{g,T})$ is continuous if and only if it has the form
\[
L(x)=cx+\sum_{j=1}^n d_jf(b_jx),\qquad c,d_j\in\R,\quad b_j\in B,
\]
and, when $c\ne0$, satisfies $\mathcal F(g,cT)\subseteq\mathcal F(f,S)$. If $c=0$ there is no additional condition.
\end{theorem}
\begin{proof}
The usual target coordinate gives the displayed form by Lemma~\ref{tauiso:maps}. For each $a\in T$, rational linearity gives
\[
g_a(Lx)=g_{ca}(x)+\sum_{j=1}^n g(ad_j)f_{b_j}(x).
\]
The sum is a rational combination of source observations. By the rational-valued assertion of the lemma, this target observation is continuous exactly when $g_{ca}\in\mathcal F(f,S)$. For $c=0$ that term is zero. The initial-topology definition proves necessity and sufficiency.
\end{proof}

\begin{theorem}[Exact isomorphism criterion and inverse]\label{tauiso:iso}
Write a map of the preceding form, with $c\ne0$, as
\[
L=m_c\circ(\mathrm{id}+\Lambda),\qquad
\Lambda(x)=\sum_{j=1}^n\sigma_j f(b_jx),\qquad \sigma_j=d_j/c.
\]
Put $W=\operatorname{span}_{\Q}\{\sigma_1,\ldots,\sigma_n\}$ and $A=(\mathrm{id}+\Lambda)|_W$. Then $L$ is a topological group isomorphism if and only if $A:W\to W$ is invertible and
\[
\mathcal F(g,cT)=\mathcal F(f,S).
\]
Its inverse is
\[
L^{-1}(y)=\frac yc-A^{-1}\Lambda(y/c)
=\frac yc-\sum_{j=1}^n(A^{-1}\sigma_j)f(b_jy/c).
\]
No additional continuity condition on the finite-rank part of this inverse is needed.
\end{theorem}
\begin{proof}
A finite sum of the Hamel observations has countable image, so an injective map on $\R$ cannot have $c=0$. The space $W$ is finite-dimensional over $\Q$, and $\Lambda(\R)\subseteq W$. On each coset $x+W$, the map $M=\mathrm{id}+\Lambda$ has the form $x+w\mapsto x+\Lambda(x)+A(w)$. Hence $M$ is bijective exactly when $A$ is invertible. Solving gives $M^{-1}(z)=z-A^{-1}\Lambda(z)$ and the displayed inverse of $L=m_c\circ M$. Division by $c$ is inside the argument of $\Lambda$; rational linearity does not allow it to be moved outside for an arbitrary real $c$.

If $L$ is an isomorphism, Theorem~\ref{tauiso:hom} gives the forward inclusion. The inverse formula has ordinary scalar coefficient $1/c$ and a countable-image remainder. Any representation of the continuous inverse supplied by Lemma~\ref{tauiso:maps} has the same ordinary coefficient: the difference of two countable-image maps has countable image, whereas a nonzero multiple of the identity has uncountable image. Applying Theorem~\ref{tauiso:hom} to the inverse yields $\mathcal F(f,S/c)\subseteq\mathcal F(g,T)$. Composition with $m_c$ gives the reverse inclusion and therefore equality.

Conversely the equality gives continuity of $L$. For each $b_j$ there is $\beta_j\in T$ with $f_{b_j}=g_{c\beta_j}$. Consequently $f(b_jy/c)=g(\beta_jy)$, so the inverse's finite-rank part is a real combination of target observations. The same equality also gives $\mathcal F(f,S/c)=\mathcal F(g,T)$. Theorem~\ref{tauiso:hom} proves continuity of the inverse.
\end{proof}

\begin{corollary}[Isomorphism types and the continuous scalar ring]\label{tauiso:types}
The spaces $(\R,\tau_{f,S})$ and $(\R,\tau_{g,T})$ are additively topologically isomorphic exactly when $\mathcal F(g,cT)=\mathcal F(f,S)$ for some real $c\ne0$. Whenever an isomorphism exists, a scalar isomorphism exists. Their continuous scalar rings therefore agree. The family $\{\tau_{f,\Q}:f\ne0\}$ has exactly $2^{\mathfrak c}$ additive topological group isomorphism types.
\end{corollary}
\begin{proof}
Necessity follows from Theorem~\ref{tauiso:iso}; for sufficiency take $\Lambda=0$. Scalar maps commute with every scalar multiplication, proving invariance of the continuous scalar ring. For rational parameter spaces the equality criterion is $g(cx)=qf(x)$ with $q\in\Q\setminus\{0\}$. There are $2^{\mathfrak c}$ nonzero rational-linear functionals, and for fixed $f$ this relation allows at most $\mathfrak c$ choices of $g$. Thus there are exactly $2^{\mathfrak c}$ classes. This count concerns the rational-parameter family; $S=\{0\}$ gives the single usual topology.
\end{proof}

\begin{corollary}[Automorphisms and shears]\label{tauiso:auto}
The automorphisms of $(\R,\tau_{f,\Q})$ are exactly
\[
L_{c,s}(x)=cx+s f(x),\qquad c\in\Q\setminus\{0\},\quad c+f(s)\ne0.
\]
Their composition and inverse are
\[
L_{c,s}\circ L_{c',s'}=L_{cc',\,cs'+c's+s f(s')},
\]
\[
L_{c,s}^{-1}(y)=\frac yc-\frac{s}{c(c+f(s))}f(y).
\]
The subgroup of shears $M_s=L_{1,s}$ is a split extension of $\Q\setminus\{0\}$ under multiplication by $(\ker f,+)$ and is nonabelian. More generally, if $1\in S$, every $M_s$ with $1+f(s)\ne0$ is an automorphism of $(\R,\tau_{f,S})$.
\end{corollary}
\begin{proof}
For $S=T=\Q$, Lemma~\ref{tauiso:maps} gives $L=c\,\mathrm{id}+sf$. Equality of the rational observation spaces forces $f_c=qf=f_q$ for some nonzero rational $q$. Injectivity of $a\mapsto f_a$ gives $c=q$. On $W=\operatorname{span}_{\Q}\{s/c\}$ the map $A$ multiplies by $1+f(s/c)$, so invertibility is equivalent to $c+f(s)\ne0$. Direct substitution gives the composition and inverse formulas.

The map $M_s\mapsto1+f(s)$ is a surjective homomorphism onto the nonzero rationals under multiplication; its kernel is $\{M_k:k\in\ker f\}$ with the additive composition law. Choose $u$ with $f(u)=1$. The assignment $q\mapsto M_{(q-1)u}$ is a splitting. For nonzero $k\in\ker f$, one has $M_uM_k=M_{u+k}$ but $M_kM_u=M_{u+2k}$, proving noncommutativity. Finally, for any $S$ containing $1$, the finite-rank term $sf$ is an allowed observation term and the scalar coefficient is $1$. Theorem~\ref{tauiso:iso} proves the stated general shear assertion.
\end{proof}

The isomorphism classification above is independent of the two-chain counterexample. No historical priority claim is made.

\section{Non-nuclearity of the prime and quotient-chain topologies}\label{sec:nuclearity}
We use nuclearity in the standard sense for abelian topological groups. Banaszczyk and Mart\'in-Peinador prove that every nuclear group respects compactness: its compact subsets are precisely those compact in its weak character topology \cite{BMP1996}. The following consequence supplies the required obstruction.

\begin{proposition}[A countable character group obstruction]\label{nuclear:countable}
Let $(G,\rho)$ be a Hausdorff abelian group whose full character group $E$ is countable and separates points. If $(G,\rho)$ is nuclear, then $\rho=\sigma(G,E)$. In particular, it is precompact.
\end{proposition}
\begin{proof}
Put $w=\sigma(G,E)$. The character embedding into $\T^E$ makes $w$ Hausdorff and metrizable, and $w\subseteq\rho$. Suppose $x_n\to x$ in $w$. The set $K=\{x\}\cup\{x_n:n\ge1\}$ is $w$-compact, hence $\rho$-compact by the cited theorem. The identity from $(K,\rho|_K)$ to $(K,w|_K)$ is a continuous bijection from a compact space to a Hausdorff space, so is a homeomorphism. Thus $x_n\to x$ in $\rho$. The identity from the metrizable space $(G,w)$ to $(G,\rho)$ is sequentially continuous, hence continuous. This proves equality. As a subgroup of $\T^E$ with its subspace topology, $(G,w)$ is precompact.
\end{proof}

To distinguish the prime from the quotient map $p$, write $\ell\in\{2,3\}$. Recall the topologies of \cite[Section 22]{HamelCompanion}:
\[
P_{\ell,d}(q)=\sup_{n\ge0}\left\|p\left(\frac{q}{d\ell^{n^2}}\right)\right\|_{\T},\qquad d\ge1,
\]
with $\gamma_\ell$ generated by the balls $P_{\ell,d}<\varepsilon$. The rational chain criterion, Theorem~\ref{chains:F}, applied to $M_n=\ell^{n^2}$, gives their full dual $\Gamma$, since consecutive ratios tend to infinity.

\begin{theorem}[Both prime topologies are non-nuclear]\label{nuclear:primes}
Neither $(\Q,\gamma_2)$ nor $(\Q,\gamma_3)$ is nuclear. More explicitly, for each $\ell\in\{2,3\}$ the sequence
\[
a_n^{(\ell)}=\frac{n!}{\ell^{v_\ell(n!)}}\,\ell^{n^2-1},\qquad n\ge1,
\]
where $v_\ell(n!)$ denotes the exponent of $\ell$ in $n!$, converges to zero in $\sigma(\Q,\Gamma)$ but not in $\gamma_\ell$.
\end{theorem}
\begin{proof}
For a prime other than $\ell$, its exponent in $a_n^{(\ell)}$ is its exponent in $n!$, which tends to infinity. The exponent of $\ell$ is $n^2-1$. Hence every fixed positive integer divides $a_n^{(\ell)}$ eventually. For every $r\in\Q$, $p(ra_n^{(\ell)})=0$ eventually, proving weak convergence. On the other hand,
\[
P_{\ell,1}(a_n^{(\ell)})\ge
\left\|p\left(\frac{n!/\ell^{v_\ell(n!)}}{\ell}\right)\right\|_{\T}
=\frac1\ell\qquad(\ell=2,3).
\]
The numerator is an integer not divisible by $\ell$, giving the stated equality for these two primes. Thus convergence fails in $\gamma_\ell$. Its full dual $\Gamma$ is countable and separates points, so Proposition~\ref{nuclear:countable} excludes nuclearity.
\end{proof}

\begin{corollary}[Both quotient-chain members are non-nuclear]\label{nuclear:chains}
The two compatible topologies $\lambda_{F_1}$ and $\lambda_{F_2}$ in Theorem~\ref{twochains:H} are both non-nuclear.
\end{corollary}
\begin{proof}
Their full dual $\Z\iota$ is countable and separates points. For each $i$, the infinite family $\{m\iota:m\in F_i\}$ consists of distinct characters and is equicontinuous, directly from the defining supremum $P_{F_i}$. A precompact abelian group has only finite equicontinuous character sets. Indeed characters extend to its compact completion, and equicontinuity gives an equicontinuous family there; by the Arzelà--Ascoli theorem it has compact closure in the discrete character group of that compact group and is therefore finite. Hence neither $\lambda_{F_i}$ is precompact. Proposition~\ref{nuclear:countable} applies.
\end{proof}

These conclusions identify both members of the specified pairs. They do not assert that every incompatible compatible pair has two non-nuclear members, and they do not decide compatibility of $\gamma_2\vee\gamma_3$.

\section{Sequences and dependence on the coefficient functional}\label{sec:sequences}
Let $\mathcal D$ be the set of nonzero rational-linear maps $f:\R\to\Q$. Here $\tau_f=\tau_{f,\R}$, whereas $\sigma_f$ is generated by the usual topology and the circle-valued observations $p\circ f_a$, $a\in\R$. The topology $\beta_f$ is the weak topology generated by $H_f$. We use the result of \cite[Section 14]{HamelCompanion} that $\beta_f$ and $\sigma_f$ have the same convergent sequences. Thus
\[
x_n\longrightarrow0\text{ in }\sigma_f
\quad\Longleftrightarrow\quad
x_n\longrightarrow0\text{ usually and }
\|p(f(ax_n))\|_{\T}\longrightarrow0\quad(a\in\R).
\]
For $\tau_f$, the corresponding observation condition is usual convergence $f(ax_n)\to0$ for every $a$.

\begin{lemma}[Finite rational spans]\label{seq:finite}
Let $b_1,\ldots,b_m$ be rationally independent and $x_n=\sum_{j=1}^m q_{n,j}b_j$, with $q_{n,j}\in\Q$. For each $f\in\mathcal D$, convergence to zero in $\tau_f$ is equivalent to $q_{n,j}\to0$ usually for every $j$. Convergence to zero in $\sigma_f$, equivalently in $\beta_f$, is equivalent to usual convergence of $x_n$ and
\[
\left\|p\left(\frac{q_{n,j}}d\right)\right\|_{\T}\longrightarrow0
\qquad(1\le j\le m,\ d\ge1).
\]
\end{lemma}
\begin{proof}
This is the finite-span result of \cite{HamelCompanion}; we recall its argument. The rational-linear map $a\mapsto(f(ab_1),\ldots,f(ab_m))$ is onto $\Q^m$. Otherwise a nonzero rational vector $(r_j)$ annihilates its image, giving $f(a\sum_j r_jb_j)=0$ for all $a$, contrary to independence and $f\ne0$. Choosing coordinate vectors, and coordinate vectors divided by $d$, proves necessity of the respective conditions. Conversely each observation is the finite sum $f(ax_n)=\sum_j q_{n,j}f(ab_j)$ with fixed rational coefficients. Coordinate convergence gives real convergence; the modular conditions give convergence modulo $\Z$ after choosing a common denominator for those coefficients. The stated usual convergence supplies the remaining coordinate for $\sigma_f$.
\end{proof}

\begin{theorem}[Sequences converging for every functional]\label{seq:universal}
For a sequence $(x_n)$ in $\R$, convergence to zero in $\tau_f$ for every $f\in\mathcal D$ holds if and only if its rational span is finite-dimensional and its coordinates in a rational basis of that span tend to zero usually.

Convergence to zero in $\sigma_f$ for every $f\in\mathcal D$, equivalently in $\beta_f$ for every such $f$, holds if and only if its rational span is finite-dimensional and, in a rational basis, the usual and modular conditions of Lemma~\ref{seq:finite} hold. The zero-dimensional span is allowed, with empty coordinate conditions.
\end{theorem}
\begin{proof}
If the rational span is infinite-dimensional, choose recursively a rationally independent subsequence $(x_{n_k})$. Such a choice can always be continued past any finite prefix, since that prefix has finite-dimensional span. Prescribe $f(x_{n_k})=1/2$ for every $k$ and extend to a rational-linear functional on $\R$ using a rational basis. It is nonzero. The observation with $a=1$ prevents convergence to zero both in the real-valued topology and in the modular topology. Equality of the modular convergent sequences excludes $\beta_f$-convergence as well. Thus convergence for every functional forces finite-dimensional span. Lemma~\ref{seq:finite} proves the remaining necessity and sufficiency.
\end{proof}

\begin{theorem}[Infinite-dimensional convergence depends on the functional]\label{seq:dependence}
There exist $f,g\in\mathcal D$ and a rationally independent sequence $(x_n)$ such that $x_n\to0$ in $\tau_f$ (and hence in $\sigma_f$ and $\beta_f$), but $x_n$ does not converge to zero in any of $\tau_g,\sigma_g,\beta_g$.
\end{theorem}
\begin{proof}
Choose a transcendental real $t$ with $0<t<1$ and put $K=\Q(t)$. Transcendence identifies $K$ with the rational-function field $\Q(T)$, which embeds in the formal Laurent-series field $\Q((T))$ by expansion at $T=0$. Every such Laurent series has only finitely many negative powers. Let $\ell:K\to\Q$ extract its constant coefficient; this is rational-linear and $\ell(1)=1$.

Choose a basis of $\R$ over $K$ containing $1$, and let $\Pi:\R\to K$ be the associated $K$-linear projection onto the coordinate of $1$. Set $f=\ell\circ\Pi$, so $f\in\mathcal D$. For fixed $a\in\R$, $K$-linearity gives
\[
f(at^n)=\ell(t^n\Pi(a)).
\]
This is the coefficient of $T^{-n}$ in the Laurent expansion of $\Pi(a)$, and is therefore zero for all sufficiently large $n$. The threshold may depend on $a$, as allowed by initial-topology convergence. Since $t^n\to0$ usually, it follows that $t^n\to0$ in $\tau_f$.

The sequence $x_n=t^n$, $n\ge1$, is rationally independent by transcendence. Prescribe $g(t^n)=1/2$ for every $n$ and extend to a rational-linear functional on $\R$. Its observation at $a=1$ is constantly $1/2$, both as a real value and modulo $\Z$. This excludes convergence to zero in $\tau_g$ and $\sigma_g$, and hence also in $\beta_g$.
\end{proof}

Thus the functional independence in the finite-span lemma does not extend to arbitrary infinite-dimensional spans. Theorem~\ref{seq:universal} classifies convergence with a universal quantifier over the functional; Theorem~\ref{seq:dependence} concerns particular functionals and does not classify convergence for an arbitrary fixed one.

\section{Size of the mixed-prime character group}\label{sec:mixed-dimension}

The remaining mixed-prime problem admits a quantitative restriction on its possible character parameters. Throughout this section, distances and dimensions on $\T$ refer to its usual circle metric, and logarithms are natural. Put

\[
P(k)=\max\{P_{2,1}(k),P_{3,1}(k)\},\qquad B_\delta=\{k\in\Z:P(k)<\delta\},
\]

\[
H=\{\alpha\in\T:k\mapsto k\alpha\text{ is continuous on }(\Z,P)\},
\]

\[
K_\varepsilon=B_\varepsilon^\triangleright=\{\alpha\in\T:\|k\alpha\|_{\T}\le 1/4\text{ for every }k\in B_\varepsilon\}.
\]

For a compact set $K\subseteq\T$, let $N(K,r)$ be the least number of circle balls of radius $r$ covering $K$. Its upper box dimension is the limit superior of $\log N(K,r)/\log(1/r)$ as $r$ decreases to zero. We use Hausdorff dimension in the same circle metric.

\begin{theorem}[Dimension zero of the mixed-prime dual]\label{mixed:dimension}
For every $0<\varepsilon<1/4$, the compact polar $K_\varepsilon$ has upper box dimension zero. More quantitatively, for integers $N\ge2$,

\[
N(K_\varepsilon,1/(50N))\le \exp\bigl(C_\varepsilon\sqrt{\log N}\bigr)
\]

for some constant $C_\varepsilon>0$ independent of $N$. Consequently $H$ has Hausdorff dimension zero.

\end{theorem}
\begin{proof}
We first count integers in small metric balls. Fix $0<\delta<1$ and an integer $L>1/\delta$. Only the observations indexed by

\[
I_N(\delta)=\{(\ell,n):\ell\in\{2,3\},\ n\ge0,\ \ell^{n^2}\le N/\delta\}
\]

need be retained for integers between zero and $N$. Their number is

\[
d_N(\delta)=2+\left\lfloor\sqrt{\frac{\log(N/\delta)}{\log2}}\right\rfloor+\left\lfloor\sqrt{\frac{\log(N/\delta)}{\log3}}\right\rfloor=O_\delta(\sqrt{\log N}).
\]

Map $0,\ldots,N$ to $\T^{d_N(\delta)}$ using the retained observations, and partition each coordinate, represented in $[0,1)$, into $L$ half-open intervals of length $1/L$. One product box contains at least $(N+1)/L^{d_N(\delta)}$ integers. Subtract its smallest integer from every integer in that box. The resulting distinct differences lie in $[0,N]$. Their retained observations have circle norm less than $1/L<\delta$, and their omitted observations have norm at most $N/\ell^{n^2}<\delta$. Thus

\[
|B_\delta\cap[0,N]|\ge\frac{N+1}{L^{d_N(\delta)}}.
\]

Now fix $\varepsilon$ as in the statement, set $\delta=\varepsilon/2$, and write $A_N=B_\delta\cap[0,N]$. If $\alpha\in K_\varepsilon$ and $k\in A_N$, subadditivity gives $P(2k)\le2P(k)<\varepsilon$. Both $k$ and $2k$ therefore belong to $B_\varepsilon$. The small-multiples lemma with two multiples yields $\|k\alpha\|_{\T}\le1/8$. Hence the trigonometric polynomial

\[
S_N(t)=\sum_{k\in A_N}e^{2\pi i kt}
\]

satisfies

\[
\operatorname{Re}S_N(\alpha)\ge |A_N|/\sqrt2\quad(\alpha\in K_\varepsilon),\qquad |S_N^{\prime}(t)|\le2\pi N|A_N|.
\]

Put $r_N=1/(100N)$. Throughout the circle interval of radius $r_N$ about each $\alpha\in K_\varepsilon$, the derivative estimate implies

\[
|S_N(t)|\ge\left(\frac1{\sqrt2}-\frac{2\pi}{100}\right)|A_N|\ge |A_N|/2.
\]

The set $K_\varepsilon$ is compact, being an intersection of closed character preimages in $\T$. Choose a maximal finite collection of its points with pairwise distances greater than $2r_N$, and denote its cardinality by $J$. The intervals of radius $r_N$ about these points are disjoint, while those of radius $2r_N$ cover $K_\varepsilon$. Orthogonality of the distinct integer exponentials gives

\[
|A_N|=\int_0^1|S_N(t)|^2\,dt\ge J(2r_N)\frac{|A_N|^2}{4}.
\]

Since $0\in A_N$, its cardinality is positive. It follows that

\[
N(K_\varepsilon,2r_N)\le J\le\frac{200N}{|A_N|}\le200L^{d_N(\delta)}.
\]

This proves the quantitative covering bound after increasing $C_\varepsilon$. Dividing its logarithm by $\log(50N)$ gives a limit of zero. Monotonicity of covering numbers between successive radii $1/(50N)$ proves upper box dimension zero. Hausdorff dimension is bounded above by upper box dimension, so each $K_\varepsilon$ also has Hausdorff dimension zero.

Finally the quarter-arc criterion gives

\[
H=\bigcup_{m\ge5}K_{1/m}.
\]

Indeed continuity supplies such a ball, and conversely a quarter-arc bound on one ball implies continuity. Countable stability of Hausdorff dimension proves the conclusion for $H$.

\end{proof}

All dyadic and triadic rational parameters belong to $H$. For example, if $n^2\ge a$, then

\[
\|k/2^a+\Z\|_{\T}\le2^{n^2-a}P(k),
\]

and the analogous estimate holds for powers of three. Since $H$ is a group, Bezout identities then give

\[
\Z[1/6]/\Z\subseteq H,\qquad \dim_{\mathrm H}H=0.
\]

In particular, almost every circle parameter, with respect to normalized Lebesgue measure, defines a discontinuous character of $(\Z,P)$. This does not prove that $H$ consists only of rational parameters: sets of Hausdorff dimension zero may contain irrational points and may be uncountable. The compatibility question still asks whether $H\subseteq\Q/\Z$. The argument above does not require a mixed-prime equidistribution estimate.

\section*{Formal verification}
The universal-sequence and functional-dependence results of Section~\ref{sec:sequences} and the dimension result of Section~\ref{sec:mixed-dimension} have Lean and Mathlib formalizations, as do the chain criteria, the two-chain construction with its transfers, and the real-valued isomorphism classification of Sections~\ref{sec:chains-audited}--\ref{sec:tauiso}. The nuclearity conclusions of Section~\ref{sec:nuclearity} rely on the cited theory of nuclear groups and remain outside the verified portion, except for the explicit weakly null sequences. Formal verification establishes the stated conclusions under their stated hypotheses; it does not certify every printed proof.

The Lean and Mathlib sources for the selected formalized results are available in the public repository \href{https://github.com/irinayar/MyPseudotopology-Lean}{MyPseudotopology-Lean}. The source snapshot cited here is \href{https://github.com/irinayar/MyPseudotopology-Lean/tree/c0444e290d5b7f0952907356f7f6a23e62b87586}{commit \texttt{c0444e2}}. Its verification records state the checked declarations and their source hashes. Individual module checks are distinguished from a full rebuild of the entire library.

\section*{Remaining questions}
The two-chain construction resolves Conjecture 23.2 and Question B.4 of \cite{HamelCompanion}. Theorem~\ref{nuclear:primes} proves that both prime topologies $\gamma_2$ and $\gamma_3$ are non-nuclear, and Corollary~\ref{nuclear:chains} identifies both quotient-chain members as non-nuclear. These are applications of the established Glicksberg theorem for nuclear groups, not a new general compactness theorem.

The full character group of the particular mixed-prime join $\gamma_2\vee\gamma_3$ remains undetermined. Equivalently, with $P(k)=\max\{P_{2,1}(k),P_{3,1}(k)\}$ on $\Z$, it remains to determine whether every continuous character of this metric group has rational parameter modulo $\Z$; see \cite[Section 22]{HamelCompanion}. Theorem~\ref{mixed:dimension} proves that the integer character parameter group has Hausdorff dimension zero, but does not exclude irrational parameters. Non-nuclearity alone does not answer this question.

For sequences, Theorem~\ref{seq:universal} classifies convergence for every coefficient functional, and Theorem~\ref{seq:dependence} shows functional dependence in infinite dimension. A classification of all infinite-dimensional rational-span sequences converging for an arbitrary fixed $f$ remains open in these arguments.

\part{Rational Scalar Continuity and Hamel Character Topologies on the Real Line}\label{part:II}
\setcounter{section}{0}

\section{Introduction and relation to earlier work}
A real vector space can carry a group topology compatible with some scalar multiplications but not others. In a real topological vector space the scalar multiplication is jointly continuous. A commonly used weaker convention for a real topological vector group requires only continuity of \(x\mapsto ax\) for each fixed real scalar \(a\). We weaken the latter condition by requiring continuity only for rational scalars. The algebraic real vector space structure remains part of the object.

The fixed-operator framework is established: Bourbaki treats topological groups with operators in \cite[Chapter III, Section 6, no. 1]{BourbakiGT}. Taking the operators to be rational scalar multiplications places our continuity requirement in that framework, equivalently in the theory of vector spaces over the discrete field \(\Q_d\). Salas and Tapia-Garc\'ia \cite{ST2016} study the related, more specialized theory of extended topological vector spaces and extended seminorms; their work is not being cited here as an origin of the general operator-group framework. Real topological vector groups and character lifting are discussed by Dom\'inguez and Tarieladze \cite{DT2008}. Some literature imposes additional local balancedness assumptions, so we state our convention explicitly.

The companion article \cite{HamelCompanion} uses all real rescalings \(f_a(x)=f(ax)\) of a rational-valued additive function. It separates modular observations from real-valued observations and identifies compact completion models. We retain the same method but permit an arbitrary nonzero rational parameter space \(S\). This permits exact control of the continuous scalar ring, including examples in which a continuous nonzero scalar map has a discontinuous inverse.

The primary source for the divisibility-chain criteria, the two-chain counterexample and its transfers, and the real-valued homomorphism and isomorphism classification is Part I (Sections~\ref{sec:chains-audited}--\ref{sec:tauiso}), which treats both the full real-orbit and the rational-parameter cases. We use those results by reference to Part I; their proofs are not repeated in this part. Sections~\ref{q:sec:chains-audited}--\ref{q:sec:tauiso} recall the statements needed for the present applications. Our rational-scalar framework is developed in the preceding sections: arbitrary rational parameter spaces, exact continuous scalar rings, sequence-separation quantifiers, rational scalar saturation, and the transfer theorem of Section~\ref{q:sec:transfer}. In particular, the TVGQ nonexistence conclusion uses the counterexample from Part I together with this rational-scalar framework.

The rational-scalar results are proved below from finite rational linear algebra, density of kernels, and standard character duality. Sections~\ref{q:sec:sequences}--\ref{q:sec:transfer} address two further issues. Universal sequence separation fails for the rational parameter space even though separation holds for each fixed functional. Compatible rational saturation identifies the greatest-compatible TVGQ question with the unrestricted Mackey problem for the same dual, and a rational transfer construction produces a strict compatible refinement without settling greatestness. These statements do not establish historical novelty. In particular, no claim is made that existing operator-group literature lacks these results.

\section{TVGQ and the continuous scalar ring}
\begin{definition}
A \emph{TVGQ} is a real vector space \(V\) with a Hausdorff topology \(\tau\) such that \((V,+,\tau)\) is a topological group and every map
\[
m_q:V\longrightarrow V,\qquad m_q(x)=qx,\qquad q\in\Q,
\]
is continuous. A morphism of TVGQ spaces will mean a continuous real-linear map. A \emph{real TVG} in this article requires continuity of every fixed real scalar map, without joint continuity in the scalar variable.
\end{definition}

The chosen morphisms retain the algebraic real action. Merely forgetting that action produces a topological vector space over \(\Q_d\), but forgets part of the object.

\begin{proposition}\label{q:prop:division}
For an additive group topology on a real vector space, the TVGQ condition is equivalent to continuity of \(x\mapsto x/n\) for every positive integer \(n\). It is also equivalent to each \(m_n\), \(n\ge1\), being a homeomorphism.
\end{proposition}
\begin{proof}
Integer multiplication is continuous in every topological abelian group. If division by \(n\) is continuous, composition with integer multiplication gives continuity of every \(m_{k/n}\). The map \(m_n\) is an algebraic bijection with inverse \(m_{1/n}\), proving the remaining equivalence.
\end{proof}

\begin{definition}
For a TVGQ space define
\[
C_\tau=\{a\in\R:m_a\text{ is continuous}\}.
\]
We call \(C_\tau\) its \emph{continuous scalar ring}.
\end{definition}

\begin{proposition}\label{q:prop:ring}
The set \(C_\tau\) is a unital subring of \(\R\) containing \(\Q\). A nonzero scalar \(a\) acts by a homeomorphism exactly when \(a,a^{-1}\in C_\tau\).
\end{proposition}
\begin{proof}
Use \(m_{a+b}=m_a+m_b\), \(m_{ab}=m_a\circ m_b\), and \(m_{-a}=-m_a\). The homeomorphism assertion follows because the algebraic inverse of \(m_a\) is \(m_{a^{-1}}\).
\end{proof}

Continuity of a bijective scalar map does not by itself imply continuity of its inverse. Consequently \(C_\tau\) should not be called a field without an additional argument.

\section{Hamel observations and finite-coordinate lemmas}
Write \(p:\R\to\T=\R/\Z\) for the quotient homomorphism, where \(\T\) has its usual compact topology, and put
\[
\|t\|_{\T}=\inf_{k\in\Z}|t-k|.
\]
Fix a nonzero rational-linear map \(f:\R\to\Q\). Such a map can be obtained by selecting one coefficient of a Hamel basis of \(\R\) over \(\Q\). Conversely any nonzero such map can be represented as a coefficient functional after choosing \(u\) with \(f(u)=1\) and a basis of \(\ker f\). It is onto \(\Q\).

For \(a\in\R\), put \(f_a(x)=f(ax)\). Then
\[
f_{a+b}=f_a+f_b,\qquad f_{qa}=qf_a\quad(q\in\Q).
\]
The map \(a\mapsto f_a\) is injective: if \(a\ne0\), multiplication by \(a\) is onto \(\R\), so \(f_a\) cannot vanish identically.

\begin{lemma}\label{q:lem:finite}
Let \(a_1,\ldots,a_n\) be rationally independent real numbers. The map
\[
F:\R\longrightarrow\Q^n,\qquad F(x)=(f_{a_1}(x),\ldots,f_{a_n}(x))
\]
is onto, and its kernel is dense in the usual real line. More generally the common kernel of finitely many rational-valued rational-linear maps on \(\R\) is usual-dense.
\end{lemma}
\begin{proof}
If \(F\) were not onto, a nonzero rational linear functional on \(\Q^n\) would annihilate its image. Thus \(f_{\sum q_ja_j}=0\) for some nonzero rational coefficient vector, contradicting injectivity and independence.

For the second assertion a common kernel has finite codimension over \(\Q\). It is nonzero because \(\R\) has infinite rational dimension. If \(k\ne0\) lies in it, it contains \(k\Q\), which is usual-dense. This also proves the density assertion for \(F\).
\end{proof}

\begin{lemma}\label{q:lem:mixed}
If \(c,a\in\R\) satisfy \(p(cx+f(ax))=0\) for every \(x\in\R\), then \(c=a=0\).
\end{lemma}
\begin{proof}
An additive map \(h:\R\to\Z\) is zero: \(h(x)=nh(x/n)\), so \(h(x)\) is divisible by every positive integer. Hence \(cx+f(ax)=0\). Its rational-valued second term forces \(c=0\), since otherwise \(c\R=\R\) would be contained in \(\Q\). Injectivity of \(a\mapsto f_a\) now gives \(a=0\).
\end{proof}

\section{Three topologies from a rational parameter space}\label{q:sec:three}
Let \(S\subseteq\R\) be a nonzero rational vector subspace. Define
\[
\begin{aligned}
\beta_{f,S}
 &=\init\{x\mapsto p(bx),\,x\mapsto p(f(ax)):b\in\R,\ a\in S\},\\
\sigma_{f,S}&=\tu\vee\beta_{f,S},\\
\tau_{f,S}
 &=\init\{\mathrm{id}:\R\to\R_{\mathrm{usual}},\
 f_a:\R\to\Q_{\mathrm{usual}}:a\in S\}.
\end{aligned}
\]
All initial topologies here use the displayed target topologies. Since \(f_a\) is rational-valued, using \(\R_{\mathrm{usual}}\) instead of \(\Q_{\mathrm{usual}}\) in the last line gives the same initial topology.

For example, basic zero-neighborhoods for \(\sigma_{f,S}\) and \(\tau_{f,S}\) respectively are
\[
\begin{aligned}
V(\delta;A,\varepsilon)&=\{x:|x|<\delta,\ \|f(ax)\|_{\T}<\varepsilon\ (a\in A)\},\\
U(\delta;A,\varepsilon)&=\{x:|x|<\delta,\ |f(ax)|<\varepsilon\ (a\in A)\},
\end{aligned}
\]
where \(A\subseteq S\) is finite and \(\delta,\varepsilon>0\).

\begin{proposition}\label{q:prop:basic}
These are Hausdorff nondiscrete TVGQ topologies. The topology \(\beta_{f,S}\) is precompact, \(\sigma_{f,S}=\tu\vee\beta_{f,S}\), and
\[
\beta_{f,S}\subseteq\sigma_{f,S}\subseteq\tau_{f,S}.
\]
\end{proposition}
\begin{proof}
All observations are homomorphisms, so their initial topologies are group topologies. Ordinary characters separate points for \(\beta_{f,S}\); the other two contain the usual topology. Rational dilation preserves the observation family because \(qS\subseteq S\), and the usual coordinate remains continuous. Precompactness follows from the embedding into a product of circles.

For each fixed \(a\), \(f(a/n)=f(a)/n\to0\). Thus \(1/n\to0\) in \(\tau_{f,S}\), and hence in the two coarser topologies, proving nondiscreteness. The stated inclusions follow directly from the definitions.
\end{proof}

\section{Completions and the full modular character groups}
Let
\[
K_{\R}=\widehat{\R_d},\qquad K_S=\widehat{S_d}.
\]
The hats denote Pontryagin duals; the subscripts \(d\) mean discrete additive groups. Both duals are compact, with their pointwise convergence topologies. Standard compact duality identifies the continuous characters of \(K_S\) with evaluations at elements of \(S\), and similarly for \(K_{\R}\); see \cite{Morris1977}.

Define
\[
e(x)(b)=p(bx),\qquad j_S(x)(a)=p(f(ax)).
\]
These give embeddings
\[
I_S(x)=(x,j_S(x)),\qquad J_S(x)=(e(x),j_S(x))
\]
for \(\sigma_{f,S}\) and \(\beta_{f,S}\), respectively.

\begin{theorem}\label{q:thm:completion}
As Hausdorff topological groups,
\[
\boxed{\Comp(\R,\sigma_{f,S})\cong\R_{\mathrm{usual}}\times K_S,\qquad
\Comp(\R,\beta_{f,S})\cong K_{\R}\times K_S.}
\]
These identifications use the embeddings \(I_S,J_S\).
\end{theorem}
\begin{proof}
A continuous character on \(\R\times K_S\) has the form
\[
(r,\kappa)\longmapsto p(cr)+\kappa(a),\qquad c\in\R,\ a\in S.
\]
On \(I_S(\R)\) it becomes \(p(cx+f(ax))\). Lemma~\ref{q:lem:mixed} says that its annihilator is zero. The same calculation applies to \(J_S(\R)\) in \(K_{\R}\times K_S\), whose characters are evaluations at pairs \((c,a)\).

A subgroup of a locally compact abelian group with zero annihilator is dense. Both targets are complete Hausdorff groups, so these dense embeddings identify the group completions.
\end{proof}

\begin{theorem}\label{q:thm:characters}
The full continuous character groups of \(\beta_{f,S}\) and \(\sigma_{f,S}\), regarded as sets of functions on \(\R\), are both
\[
\boxed{H_{f,S}=\{\chi_{c,a}:c\in\R,\ a\in S\},\qquad
\chi_{c,a}(x)=p(cx+f(ax)).}
\]
The parameters are unique. Furthermore
\[
\Hom_{\mathrm{cont}}((\R,\sigma_{f,S}),\R_{\mathrm{usual}})
 =\R\,\mathrm{id},\qquad
\Hom_{\mathrm{cont}}((\R,\beta_{f,S}),\R_{\mathrm{usual}})=\{0\}.
\]
\end{theorem}
\begin{proof}
A continuous homomorphism into a complete Hausdorff group extends uniquely to the group completion. Apply Theorem~\ref{q:thm:completion} with target \(\T\); its product character descriptions give exactly \(H_{f,S}\), and Lemma~\ref{q:lem:mixed} gives uniqueness.

For real-valued maps, every continuous homomorphism from a compact group into \(\R\) is zero, because a compact subgroup of \(\R\) is trivial. Thus the compact completion of \(\beta_{f,S}\) gives zero, while \(\R\times K_S\) contributes precisely the usual real linear maps on its real factor.
\end{proof}

The equality of character groups makes the modular pair \emph{compatible}. Their compact-open dual topologies are not identified by this assertion.

For the real-valued construction choose a rational basis \(B\) of \(S\). Rational linear combinations are continuous, so the maps \(f_b\), \(b\in B\), generate the same topology as all maps \(f_a\), \(a\in S\).

\begin{theorem}\label{q:thm:realcompletion}
The graph embedding
\[
G_B(x)=(x,(f(bx))_{b\in B})
\]
identifies
\[
\boxed{\Comp(\R,\tau_{f,S})\cong
\R_{\mathrm{usual}}\times\R_{\mathrm{usual}}^{B}.}
\]
Consequently the continuous real-valued additive maps on \((\R,\tau_{f,S})\) are exactly
\[
x\longmapsto cx+\sum_{j=1}^n d_jf(b_jx),
\qquad c,d_j\in\R,\quad b_j\in B.
\]
\end{theorem}
\begin{proof}
The graph induces exactly \(\tau_{f,S}\). To prove density, take a basic product open set, restricting only finitely many coordinates \(b_1,\ldots,b_n\) besides the first. Choose rational values in their specified real intervals. Lemma~\ref{q:lem:finite} gives a fiber with those values and a usual-dense kernel. Adjusting within that fiber places its first coordinate in the prescribed interval. Thus the graph is dense. The product of complete real groups is complete.

A continuous additive map from \(\R\times\R^B\) into \(\R\) depends on finitely many coordinates. Indeed, a basic neighborhood mapped into \((-1,1)\) leaves all other coordinates unrestricted, so the map vanishes on that entire coordinate subgroup by integer scaling. On the remaining finite real product it is a continuous additive map and hence a real linear functional. Extend to the completion and restrict to \(G_B(\R)\).
\end{proof}

\section{Strict separation of modular and real-valued observations}
\begin{lemma}\label{q:lem:usual}
If \(g:\R\to\Q\) is nonzero and rational-linear, neither \(g\) nor \(p\circ g\) is usual-continuous.
\end{lemma}
\begin{proof}
A usual-continuous additive real function is \(x\mapsto cx\), which can be rational-valued only for \(c=0\). A usual-continuous character is \(x\mapsto p(cx)\). If \(p(g(x))=p(cx)\), the additive integer-valued difference is zero, reducing to the first assertion.
\end{proof}

\begin{theorem}\label{q:thm:strict}
For every nonzero \(S\),
\[
\boxed{\tu\subsetneq\sigma_{f,S}\subsetneq\tau_{f,S}.}
\]
Every nonzero rational-valued additive map \(g\) is unbounded on every \(\sigma_{f,S}\)-neighborhood of zero. Hence the continuous rational-valued additive duals of both modular topologies are zero. The topology \(\beta_{f,S}\) is incomparable with \(\tu\).
\end{theorem}
\begin{proof}
Choose \(a\in S\setminus\{0\}\). The character \(p\circ f_a\) is not usual-continuous by Lemma~\ref{q:lem:usual}, proving the first strict inclusion.

Take a basic neighborhood \(V(\delta;A,\varepsilon)\), and put
\[
N=\ker g\cap\bigcap_{a\in A}\ker f_a.
\]
It is usual-dense by Lemma~\ref{q:lem:finite}. Choose \(y\) with \(g(y)\ne0\), and a positive integer \(d\) clearing the denominators of all \(f_a(y)\), \(a\in A\). For each positive integer \(k\), choose \(z_k\in N\) such that \(x_k=kdy+z_k\) satisfies \(|x_k|<\delta\). Then
\[
f_a(x_k)=kd f_a(y)\in\Z,\qquad
|g(x_k)|=kd|g(y)|\to\infty.
\]
Thus \(x_k\in V\) and \(g\) is unbounded there. No such \(g\) is continuous. Since \(f_a\) is continuous for \(\tau_{f,S}\), the second inclusion is strict. A map continuous for the coarser \(\beta_{f,S}\) would also be continuous for \(\sigma_{f,S}\), proving its rational-dual assertion.

The character already used shows \(\beta_{f,S}\not\subseteq\tu\). Conversely, if \(\tu\subseteq\beta_{f,S}\), the usual real group would be a continuous homomorphic image of a precompact group and therefore precompact, which it is not.
\end{proof}

\section{Exact continuous scalar rings}
For a nonzero rational subspace \(S\), write
\[
M(S)=\{r\in\R:rS\subseteq S\}.
\]
This is a unital rational subalgebra of \(\R\).

\begin{theorem}[Exact scalar realization]\label{q:thm:scalar}
For each \(\eta=\beta_{f,S},\sigma_{f,S},\tau_{f,S}\),
\[
\boxed{C_\eta=M(S).}
\]
For \(r\ne0\), multiplication by \(r\) is a homeomorphism exactly when \(rS=S\).
\end{theorem}
\begin{proof}
If \(rS\subseteq S\), precomposition by \(m_r\) replaces \(f_a\) by \(f_{ar}\); ordinary characters and the usual coordinate remain continuous. This proves sufficiency in all three constructions.

For either modular topology, suppose \(m_r\) is continuous. For each \(a\in S\), the character \(p\circ f_{ar}\) is continuous. Theorem~\ref{q:thm:characters} gives
\[
p(f(arx))=p(cx+f(bx)),\qquad c\in\R,\ b\in S.
\]
The difference is additive and integer-valued, so \(f_{ar-b}(x)=cx\). Rational-valuedness forces \(c=0\), and injectivity of the parameter map gives \(ar=b\in S\).

For \(\tau_{f,S}\), continuity of \(m_r\) makes each \(f_{ar}\) a continuous real-valued additive map. Theorem~\ref{q:thm:realcompletion} yields
\[
f_{ar}(x)=cx+\sum_{j=1}^n d_j f_{b_j}(x)
\]
with distinct basis elements \(b_j\in B\). On their common kernel, \(cx\) is rational-valued. That kernel has finite rational codimension and is uncountable: otherwise \(\R\), a countable union of its cosets, would be countable. Thus \(c=0\). Surjectivity onto \(\Q^n\) now shows each \(d_j\in\Q\), by testing the rational standard unit vectors. Injectivity gives \(ar=\sum_jd_jb_j\in S\). If the list is empty, rational-valuedness directly gives \(c=0\) and then \(ar=0\). Therefore \(rS\subseteq S\).

Finally Proposition~\ref{q:prop:ring} gives the homeomorphism criterion: \(rS\subseteq S\) and \(r^{-1}S\subseteq S\) are equivalent to \(rS=S\).
\end{proof}

\begin{corollary}\label{q:cor:algebra}
Every unital \(\Q\)-subalgebra \(A\subseteq\R\) occurs as the exact continuous scalar ring of each of the three constructions. Every field \(\Q\subseteq F\subseteq\R\) occurs in particular.
\end{corollary}
\begin{proof}
Take \(S=A\). If \(rA\subseteq A\), then \(r=r\cdot1\in A\); the reverse inclusion follows from multiplication in \(A\). Thus \(M(A)=A\).
\end{proof}

\begin{corollary}[Continuous bijection with discontinuous inverse]
Let \(t\in\R\) be transcendental and take \(S=\Q[t]\). In each construction, multiplication by \(t\) is continuous and bijective, but its inverse is discontinuous.
\end{corollary}
\begin{proof}
The continuous scalar ring is \(\Q[t]\). It contains \(t\), but not \(1/t\): a polynomial identity \(tP(t)=1\) would contradict transcendence.
\end{proof}

\section{Order of the parameter spaces and the rational case}
\begin{theorem}\label{q:thm:order}
For nonzero rational subspaces \(S,T\subseteq\R\) and a fixed \(f\), each of the following is equivalent to \(S\subseteq T\):
\[
\beta_{f,S}\subseteq\beta_{f,T},\qquad
\sigma_{f,S}\subseteq\sigma_{f,T},\qquad
\tau_{f,S}\subseteq\tau_{f,T}.
\]
\end{theorem}
\begin{proof}
Inclusion of parameter spaces gives inclusion of the initial topologies. For the converse in either modular case, every \(p\circ f_a\), \(a\in S\), belongs to \(H_{f,T}\). The no-mixed-relation argument in the proof of Theorem~\ref{q:thm:scalar} gives \(a\in T\).

For the real-valued case every \(f_a\), \(a\in S\), is continuous for \(\tau_{f,T}\). Apply the finite-coordinate formula of Theorem~\ref{q:thm:realcompletion}. The same uncountable-kernel and rational-unit-vector argument used in Theorem~\ref{q:thm:scalar} yields \(a\in T\).
\end{proof}

\begin{corollary}[The rational case]\label{q:cor:rational}
For \(S=\Q\), the topologies \(\sigma_{f,\Q}\) and \(\tau_{f,\Q}\) are metrizable, while \(\beta_{f,\Q}\) is not metrizable. Moreover
\[
\begin{gathered}
C_{\beta_{f,\Q}}=C_{\sigma_{f,\Q}}=C_{\tau_{f,\Q}}=\Q,\\
\tu\subsetneq\sigma_{f,\Q}\subsetneq\tau_{f,\Q},\\
\Comp(\R,\beta_{f,\Q})\cong K_{\R}\times K_{\Q},\quad
\Comp(\R,\sigma_{f,\Q})\cong\R\times K_{\Q},\quad
\Comp(\R,\tau_{f,\Q})\cong\R^2.
\end{gathered}
\]
\end{corollary}
\begin{proof}
The scalar and completion formulas are the preceding theorems, since a rational basis of \(\Q\) is \(\{1\}\). The compact group \(K_{\Q}\) is metrizable because its discrete dual \(\Q\) is countable, so the embeddings give metrizability of \(\sigma_{f,\Q}\) and \(\tau_{f,\Q}\).

The compact group \(K_{\R}\) is not metrizable because its discrete dual is uncountable. A metrizable precompact group has metrizable completion: an invariant compatible metric has the group uniformity and its metric completion is the group completion. The displayed completion for \(\beta_{f,\Q}\) therefore excludes metrizability.
\end{proof}

\begin{remark}
The precompact construction retains all ordinary characters, and that uncountable family produces the nonmetrizable factor \(K_{\R}\), even when the Hamel parameter space is countable.
\end{remark}

For \(F=\Q(\sqrt2)\), the continuous scalar ring is exactly \(F\); for \(F=\R\), the full real-orbit constructions of \cite{HamelCompanion} are recovered. Different intermediate parameter spaces give different topologies by Theorem~\ref{q:thm:order}. In particular the rational constructions are strictly coarser than their real-orbit counterparts.

\section{Classification of precompact TVGQ topologies}
Let \(\mathcal A=\Hom((\R,+),\T)\) be the group of all algebraic characters. It carries a rational vector space structure
\[
(q\cdot\chi)(x)=\chi(qx),\qquad q\in\Q.
\]
Scalar multiplication here acts on the argument, not by multiplication of a circle value by a rational number. Additivity of \(\chi\) verifies the vector identities.

For a subgroup \(H\subseteq\mathcal A\), write \(\theta_H\) for its initial character topology. We use the standard precompact representation: if \(H\) separates points, \(\theta_H\) is Hausdorff precompact and its full character group is exactly \(H\); every Hausdorff precompact abelian group is generated by its continuous characters. One proof embeds the group in \(\T^H\), takes its compact closure, extends a character from that compact subgroup to the product, and uses that product characters are finite integer combinations of coordinate characters. These are standard consequences of compact duality \cite{Morris1977}.

\begin{theorem}\label{q:thm:classification}
The assignment \(H\mapsto\theta_H\) is an order-preserving bijection between point-separating rational vector subspaces of \(\mathcal A\) and Hausdorff precompact TVGQ topologies on the fixed real line. Moreover
\[
C_{\theta_H}=\{a\in\R:\chi\circ m_a\in H\text{ for every }\chi\in H\}.
\]
\end{theorem}
\begin{proof}
For a point-separating rational subspace, the standard representation gives Hausdorffness, precompactness, and the exact dual \(H\). Rational invariance makes \(m_q\) continuous because precomposition preserves all generating characters.

Conversely the dual of a Hausdorff precompact TVGQ separates points, generates its topology, and is invariant under rational precomposition. The exact dual identity proves injectivity and order reflection. Finally a map into an initial character topology is continuous exactly when its compositions with all generating characters are continuous, proving the scalar formula.
\end{proof}

Unlike in the real-scalar classification, a nonzero rational subspace need not separate points. For example \(\Q\cdot(p\circ f)\) has common kernel \(\ker f\).

\begin{proposition}[A minimal rational character topology]\label{q:prop:minimal}
Put
\[
H_{\Q}=\{x\mapsto p(qx):q\in\Q\},\qquad
\theta_{\Q}=\theta_{H_{\Q}}.
\]
Then \(\theta_{\Q}\) is a Hausdorff metrizable precompact TVGQ topology, \(C_{\theta_{\Q}}=\Q\), and \(\Comp(\R,\theta_{\Q})\cong K_{\Q}\). It is minimal among Hausdorff TVGQ topologies. It is strictly coarser than the usual Bohr topology \(t=\init\{p(bx):b\in\R\}\).
\end{proposition}
\begin{proof}
For \(x\ne0\), choose a rational \(q\) with \(0<|qx|<1\), proving separation. The countable character family gives a metrizable embedding into \(\T^{\Q}\). The map
\[
x\longmapsto(q\mapsto p(qx))
\]
has dense image in \(K_{\Q}\): its annihilator consists of rational \(q\) satisfying \(p(qx)=0\) for every real \(x\), hence only zero. This gives the completion.

By Theorem~\ref{q:thm:classification}, continuity of \(m_a\) forces \(x\mapsto p(ax)\in H_{\Q}\), so \(p((a-q)x)=0\) for some rational \(q\) and every \(x\), whence \(a=q\). The converse is automatic.

Any coarser Hausdorff TVGQ topology is precompact, so corresponds to a point-separating rational subspace of the one-dimensional space \(H_{\Q}\). The only possibility is \(H_{\Q}\) itself. Finally \(H_{\Q}\) is a proper subgroup of the ordinary real-character group, so exact precompact duality gives \(\theta_{\Q}\subsetneq t\).
\end{proof}

Thus minimality in the TVG class can be lost when the class is enlarged to TVGQ. The topology \(t\), minimal in the fixed-real-scalar setting, already has the proper Hausdorff TVGQ coarsening \(\theta_{\Q}\).

\section{Sequence separation and its quantifiers}\label{q:sec:sequences}
For fixed \(f,S\), convergence of \(x_n\) to zero in \(\sigma_{f,S}\) means usual convergence together with \(\|f(ax_n)\|_{\T}\to0\) for each \(a\in S\). For \(\tau_{f,S}\), the latter condition is replaced by \(f(ax_n)\to0\) in the usual rational topology. These tests follow directly from the initial topology definitions. We distinguish a separating sequence chosen after \(f\) from a single sequence required to work for every nonzero \(f\).

Put \(\alpha=3-2\sqrt2\). For positive integers \(n\), define integers \(P_n,Q_n>0\) by
\[
(3+2\sqrt2)^{n^2}=P_n+Q_n\sqrt2,\qquad
z_n=n!\alpha^{n^2}=n!(P_n-Q_n\sqrt2).
\]
Conjugation gives the second identity. We have \(z_n\to0\) usually, since decay exponential in \(n^2\) dominates \(n!\), and \(n!Q_n\to\infty\).

\begin{theorem}[Separation for every fixed functional]\label{q:thm:scaled-sequence}
For every nonzero \(f:\R\to\Q\) and nonzero rational subspace \(S\), there is \(y\ne0\) such that \(x_n=yz_n\) converges to zero in \(\beta_{f,S}\) and \(\sigma_{f,S}\), but not in \(\tau_{f,S}\). For some \(a_0\in S\), \(|f(a_0x_n)|\to\infty\). The set \(\{0,x_1,x_2,\ldots\}\) is compact for both modular topologies and not compact for \(\tau_{f,S}\).
\end{theorem}
\begin{proof}
Choose \(a_0\in S\setminus\{0\}\). Since \(f_{a_0}\) is onto, choose \(y\) with \(f(a_0y)=1\). For every fixed \(a\in S\),
\[
f(ax_n)=n!\bigl(P_nf(ay)-Q_nf(ay\sqrt2)\bigr).
\]
The two functional values are rational. Their denominators divide \(n!\) eventually, so this expression is eventually an integer. Together with usual convergence, this proves modular convergence for \(\sigma_{f,S}\), hence also for its coarser topology \(\beta_{f,S}\).

Set \(r=f(a_0y\sqrt2)\in\Q\). The distinguished observation satisfies
\[
f(a_0x_n)=n!Q_n(\sqrt2-r)+z_n.
\]
Since \(\sqrt2-r\ne0\), its absolute value tends to infinity. Thus there is no \(\tau_{f,S}\)-convergence. A convergent sequence and its limit form a compact set. If this set were \(\tau_{f,S}\)-compact, its continuous bijective identity map onto the usual Hausdorff subspace would be a homeomorphism. Usual convergence would then imply \(\tau_{f,S}\)-convergence, a contradiction.
\end{proof}

\begin{theorem}[Exact criterion for the unscaled Pell sequence]\label{q:thm:pell-criterion}
The sequence \(z_n\) converges to zero in \(\sigma_{f,S}\) and \(\beta_{f,S}\) for every \(f\). It fails to converge in \(\tau_{f,S}\) exactly when the restriction of \(f\) to \(S+\sqrt2S\) is nonzero. Consequently it separates the modular and real-valued topologies for every nonzero \(f\) if and only if
\[
S+\sqrt2S=\R.
\]
\end{theorem}
\begin{proof}
The preceding eventual-integrality calculation works with \(y=1\), proving modular convergence. If \(f\) vanishes on \(S+\sqrt2S\), every \(f(az_n)\) is zero, so \(z_n\to0\) in \(\tau_{f,S}\). Otherwise choose \(a\in S\) such that \(u=f(a)\) and \(v=f(a\sqrt2)\) are not both zero. Then
\[
f(az_n)=u z_n+n!Q_n(\sqrt2u-v).
\]
For rational \(u,v\), \(\sqrt2u-v\ne0\), so the absolute value diverges. Finally a proper rational subspace of \(\R\) has a nonzero rational-linear functional vanishing on it, obtained by extending a basis. Thus every nonzero \(f\) is nonzero on \(S+\sqrt2S\) exactly when that subspace is all of \(\R\).
\end{proof}

\begin{theorem}[Obstruction to any universal sequence]\label{q:thm:no-universal}
Write \(\mathfrak c=|\R|\). Suppose \(0<\dim_{\Q}S<\mathfrak c\). For every usually null sequence \((x_n)\), there is a nonzero rational-linear \(f:\R\to\Q\) for which
\[
f(ax_n)=0\quad(a\in S,\ n\ge1).
\]
In particular \(x_n\to0\) in \(\tau_{f,S}\), and no sequence separates \(\sigma_{f,S}\) from \(\tau_{f,S}\) simultaneously for every nonzero \(f\). This applies to \(S=\Q\), and to every countable intermediate field.
\end{theorem}
\begin{proof}
Let \(W\) be the rational span of \(\{ax_n:a\in S,n\ge1\}\). If \(B\) is a rational basis of \(S\), the countably many sets \(Bx_n\) span \(W\). Therefore
\[
\dim_{\Q}W\le\max\{\aleph_0,\dim_{\Q}S\}<\mathfrak c
=\dim_{\Q}\R.
\]
Thus \(W\) is proper. Extend its basis to a Hamel basis of \(\R\), and take a nonzero coefficient functional supported on a basis element outside \(W\). It vanishes on all the stated observations, and usual convergence supplies the remaining coordinate for \(\tau_{f,S}\)-convergence. A universal \(\sigma_{f,S}\)-null sequence would necessarily be usually null, so this excludes every candidate.
\end{proof}

For \(S=\Q\), the conclusion is therefore a genuine change of quantifiers: every fixed \(f\) has a separating sequence, but there is no sequence separating for every \(f\). The obstruction is not specific to the Pell formula.

\begin{lemma}[Universal modular convergence forces finite rational span]\label{q:lem:finite-universal}
Suppose \(S\ne\{0\}\) and a sequence \((x_n)\) converges to zero in \(\sigma_{f,S}\) for every nonzero rational-linear \(f:\R\to\Q\). Then the rational span of its terms is finite-dimensional.
\end{lemma}
\begin{proof}
Otherwise choose a rationally independent subsequence \((x_{n_j})\), and fix \(a_0\in S\setminus\{0\}\). The elements \(a_0x_{n_j}\) are independent. A rational-linear functional can be prescribed to take value \(1/2\) on all of them and then extended to \(\R\) using a Hamel basis. It is nonzero and satisfies \(p(f(a_0x_{n_j}))=1/2+\Z\), contradicting the required modular convergence.
\end{proof}

\begin{theorem}[Complete criterion for a universal separating sequence]\label{q:thm:universal-complete}
Let \(S\subseteq\R\) be a nonzero rational subspace, of any dimension. There exists a single sequence converging to zero in \(\sigma_{f,S}\) and failing to converge to zero in \(\tau_{f,S}\) for every nonzero rational-linear \(f\) if and only if there are finitely many real numbers \(b_1,\ldots,b_m\) such that
\[
b_1S+\cdots+b_mS=\R.
\]
The sequence furnished in the sufficient direction also converges in every \(\beta_{f,S}\).
\end{theorem}
\begin{proof}
For necessity, Lemma~\ref{q:lem:finite-universal} makes the span \(W\) of the sequence finite-dimensional. If the rational span \(SW\) of all products \(sw\) were proper, a nonzero functional vanishing on it would make every observation of the sequence zero. Together with usual convergence, that would give \(\tau_{f,S}\)-convergence, a contradiction. A basis of \(W\) therefore supplies the required finite list.

For sufficiency, discard zero parameters, divide the list by a nonzero member, and take a rational basis of its span. The resulting finite list can be written \(1,\theta_2,\ldots,\theta_m\), with every \(\theta_j\) irrational. If its dimension is one, the covering condition says \(S=\R\), and the Pell sequence above applies. Otherwise
\[
S+\theta_2S+\cdots+\theta_mS=\R.
\]
For each \(j\), choose integers \(p_{j,n}\), positive integers \(q_{j,n}>n\), with
\[
0<|p_{j,n}-q_{j,n}\theta_j|<\frac{1}{n n!},\qquad
w_{j,n}=n!(p_{j,n}-q_{j,n}\theta_j).
\]
Such approximations follow from the elementary pigeonhole approximation argument for an irrational number. Indeed, comparing fractional parts of \(0,\theta,\ldots,N\theta\) gives an error at most \(1/N\). As \(N\) grows, arbitrarily small errors cannot use only a bounded set of positive denominators, because irrationality makes the minimum error for finitely many denominators positive. This proves the stated choices. We have \(w_{j,n}\to0\) usually and \(n!q_{j,n}\to\infty\).

For any fixed \(f,a,j\),
\[
f(aw_{j,n})=n!\bigl(p_{j,n}f(a)-q_{j,n}f(a\theta_j)\bigr)
\]
is eventually an integer. Finitely interleave these sequences, taking \(w_{2,n},\ldots,w_{m,n}\) at stage \(n\). Eventual integrality for a fixed observation is uniform over this finite list of \(j\)'s, so the interleaved sequence converges in every \(\sigma_{f,S}\) and \(\beta_{f,S}\).

For a nonzero \(f\), the covering condition supplies some \(a\in S\) and \(j\) such that \(u=f(a)\) and \(v=f(a\theta_j)\) are not both zero. Then
\[
f(aw_{j,n})=u w_{j,n}+n!q_{j,n}(\theta_j u-v).
\]
The nonzero rational pair \((u,v)\) and irrationality of \(\theta_j\) make the last coefficient nonzero. This observation diverges in absolute value on that subsequence, excluding convergence in \(\tau_{f,S}\).
\end{proof}

\begin{corollary}[Intermediate fields]\label{q:cor:field-universal}
For an intermediate field \(F\), a universal separating sequence exists if and only if \([\R:F]<\infty\), equivalently \(F=\R\). Thus no proper intermediate field admits a universal separating sequence. More generally, every nonzero rational subspace of finite rational codimension in \(\R\) admits such a sequence.
\end{corollary}
\begin{proof}
For a field, the finite covering condition is exactly finite-dimensionality of \(\R\) over \(F\). If this degree is finite, then \(\mathbb C\) is an algebraic closure of \(F\) of finite degree. The Artin--Schreier theorem gives \([\mathbb C:F]=2\); since \([\mathbb C:\R]=2\), the tower formula forces \(F=\R\); see \cite{ConradAS}. For a rational subspace \(S\) of finite codimension, start with \(S\) and add scalar multiples of \(S\) whenever the current sum is proper. Such a multiple always increases the sum: for \(s_0\in S\setminus\{0\}\) and \(r\) outside the sum, \((r/s_0)S\) contains \(r\). The finite codimension can decrease only finitely many times, so the procedure reaches \(\R\).
\end{proof}

The full-dimension case mentioned in the scope discussion is therefore decided by Theorem~\ref{q:thm:universal-complete}, rather than by dimension alone. Both outcomes occur even when \(S+\sqrt2S\ne\R\). For a positive example take \(K=\Q(\sqrt2)\) and a hyperplane \(S\) of \(\R\) as a \(K\)-vector space. Then \(\dim_{\Q}S=\mathfrak c\), \(\sqrt2S=S\), and \(S\) has rational codimension two, so a universal sequence exists. For a negative example choose a transcendence basis \(B\) of \(\R/\Q\), remove one element \(t\), and let \(F\) be the field of real numbers algebraic over \(\Q(B\setminus\{t\})\). It has rational dimension \(\mathfrak c\), contains \(\sqrt2\), and \(t\) is transcendental over \(F\). Thus \([\R:F]=\infty\), \(F+\sqrt2F=F\ne\R\), and no universal separating sequence exists. These examples also show why the failure of the particular unscaled Pell sequence does not imply failure of every universal sequence.

\section{Compatible topologies and rational scalar saturation}\label{q:sec:saturation}
Put \(\T_+=[-1/4,1/4]+\Z\). A subset \(U\) of an abelian topological group is quasi-convex if every point outside \(U\) is separated from it by a continuous character sending \(U\) into \(\T_+\) and that point outside \(\T_+\). A topology is locally quasi-convex if such neighborhoods form a local basis. Products and subgroups preserve local quasi-convexity. All three observation topologies above are locally quasi-convex, by their embeddings into products of locally compact abelian groups.

For a point-separating character subgroup \(H\) on a rational vector space \(G\), let \(\mathcal L(H)\) be all Hausdorff locally quasi-convex group topologies on \(G\) with full character set \(H\). Let \(\mathcal C_{\Q}(H)\) be the members admitting all fixed rational scalar maps continuously. For \(G=\R\), these are TVGQ topologies. A greatest member of \(\mathcal L(H)\), when it exists, is a Mackey group topology for the dual pair. The general problem and counterexamples are discussed in \cite{Gab2022,Gab2018}.

\begin{lemma}[Directed compatible joins]\label{q:lem:directed}
Every nonempty directed family in \(\mathcal L(H)\) has its join in \(\mathcal L(H)\). The analogous statement holds in \(\mathcal C_{\Q}(H)\).
\end{lemma}
\begin{proof}
The diagonal embedding into the product proves that the join is a Hausdorff locally quasi-convex group topology; it also preserves continuity of fixed rational maps when all factors have it. All characters in \(H\) remain continuous. If a character \(\chi\) is continuous for the join, it sends some basic neighborhood into \(\T_+\). This neighborhood involves finitely many family members, dominated by one member because the family is directed. It is thus a neighborhood for that member.

We use the usual quarter-arc criterion: an algebraic character sending one zero-neighborhood into \(\T_+\) is continuous. Indeed, intersect its neighborhood preimages under integer multiplications \(1,\ldots,N\). If the first \(N\) multiples of a circle value lie in \(\T_+\), its distance to zero is at most \(1/(4N)\). These intersections prove continuity. Applying the criterion shows \(\chi\in H\), so the join introduces no new characters.
\end{proof}

\begin{theorem}[Compatible rational scalar saturation]\label{q:thm:saturation}
Suppose \(H\) is invariant under rational rescaling of arguments. Every \(\eta\in\mathcal L(H)\) has a compatible rational scalar refinement. More precisely, let \(\eta_n\) be the pullback of \(\eta\) under \(x\mapsto \frac{x}{n!}\). Then
\[
\eta=\eta_1\subseteq\eta_2\subseteq\cdots,\qquad
\eta^{[\Q]}=\bigvee_{n\ge1}\eta_n\in\mathcal C_{\Q}(H).
\]
This is the least rational-scalar-continuous group topology containing \(\eta\).
\end{theorem}
\begin{proof}
The map \(m_{\frac{1}{n!}}\) is an algebraic bijection, so the pullback topology is topologically isomorphic to \(\eta\), and its dual is
\(\{\chi\circ m_{\frac{1}{n!}}:\chi\in H\}=H\).
It is therefore compatible and locally quasi-convex. Since
\(m_{\frac{1}{n!}}=m_{n+1}\circ m_{\frac{1}{(n+1)!}}\)
and integer multiplication is always continuous for a group topology, the displayed chain is increasing. Lemma~\ref{q:lem:directed} proves compatibility of its join.

Let \(q=\frac{a}{b}\in\Q\), with integers \(a\), \(b>0\). For each \(n\), choose \(k\) so that \(bn!\) divides \(k!\). Then
\[
m_{\frac{1}{n!}}\circ m_q
=m_{\frac{a k!}{bn!}}\circ m_{\frac{1}{k!}}.
\]
The integer map on the right is \(\eta\)-continuous, and the last map is continuous from the join into \(\eta\). Hence \(m_q\) is continuous into each \(\eta_n\), and therefore into their join. Finally any rational-scalar-continuous topology finer than \(\eta\) makes every \(m_{\frac{1}{n!}}\) continuous into \(\eta\); it contains every \(\eta_n\) and their join. This proves leastness.
\end{proof}

\begin{corollary}[Equivalence with the unrestricted Mackey question]\label{q:cor:mackey-equivalence}
Under the hypotheses of Theorem~\ref{q:thm:saturation}, \(\mathcal C_{\Q}(H)\) has a greatest member if and only if \(\mathcal L(H)\) has a greatest member. In that event the two greatest members coincide. Every maximal member of \(\mathcal L(H)\) already belongs to \(\mathcal C_{\Q}(H)\).
\end{corollary}
\begin{proof}
If \(\mu\) is greatest in the rational class, every \(\eta\in\mathcal L(H)\) satisfies \(\eta\subseteq\eta^{[\Q]}\subseteq\mu\). Thus \(\mu\) is greatest in the full class. Conversely a greatest member of the full class equals its compatible rational saturation, so lies in the rational class. A maximal member likewise equals its saturation.
\end{proof}

This saturation argument concerns rational scalars; it does not justify arbitrary real saturation. The real-scalar counterexample follows instead from the separate transfer proved in Theorem~\ref{twochains:transfer} of Part I.

\section{The greatest compatible TVGQ question}\label{q:sec:greatest}
For the present Hamel family set
\[
\mathcal C_{f,S}^{\Q}=\mathcal C_{\Q}(H_{f,S}).
\]
It is nonempty because it contains \(\beta_{f,S}\) and \(\sigma_{f,S}\). Every member contains \(\beta_{f,S}\), the topology generated by the prescribed character group. Those additionally containing \(\tu\) contain \(\sigma_{f,S}\). These are least-member assertions, not greatest-member assertions.

\begin{proposition}[Scalar restrictions throughout the compatible class]\label{q:prop:compatible-scalars}
For every \(\eta\in\mathcal C_{f,S}^{\Q}\),
\[
\Q\subseteq C_\eta\subseteq M(S).
\]
In particular every member of \(\mathcal C_{f,\Q}^{\Q}\) has continuous scalar ring exactly \(\Q\).
\end{proposition}
\begin{proof}
The left inclusion is its defining TVGQ property. If \(m_r\) is continuous, then \(p\circ f_{ar}\) is a character in \(H_{f,S}\) for every \(a\in S\). The proof of Theorem~\ref{q:thm:scalar}, using only the character formula, forces \(ar\in S\). Thus \(r\in M(S)\). For \(S=\Q\), this ring is \(\Q\).
\end{proof}

\begin{theorem}[The precise greatest-member criterion]\label{q:thm:greatest-criterion}
Every member of \(\mathcal C_{f,S}^{\Q}\) has a maximal extension in that class. The following are equivalent: the class has a greatest member; every binary join of its members remains compatible; its maximal member is unique; and \(\mathcal L(H_{f,S})\) has a Mackey topology.
\end{theorem}
\begin{proof}
Chains have compatible upper bounds by Lemma~\ref{q:lem:directed}, so Zorn's lemma gives maximal extensions. A greatest member bounds each binary join, which therefore cannot acquire characters outside \(H_{f,S}\). Conversely, if all binary joins remain compatible, all finite joins are members; their directed join over the entire class is compatible and greatest. Uniqueness of a maximal member is equivalent to greatestness because every member extends to a maximal one. The final equivalence is Corollary~\ref{q:cor:mackey-equivalence}, since \(H_{f,S}\) is rationally invariant.
\end{proof}

\begin{conjecture}[Greatest compatible TVGQ topology for the rational Hamel pair]\label{q:conj:greatest-q}
For every nonzero \(f:\R\to\Q\), the class \(\mathcal C_{f,\Q}^{\Q}\) has a greatest member.
\end{conjecture}
Theorem~\ref{twochains:transfer} of Part I disproves this positive assertion for every nonzero $f$. Its two compatible TVGQ topologies have an incompatible join, violating the criterion in Theorem~\ref{q:thm:greatest-criterion}. The corresponding unrestricted dual pair $(\R,H_{f,\Q})$ therefore has no Mackey topology.

\section{A rational transfer construction and a strict compatible refinement}\label{q:sec:transfer}
Write
\[
\Gamma=\{q\mapsto p(rq):r\in\Q\},\qquad
\gamma_0=\init\{q\mapsto p(rq):r\in\Q\}
\]
on the additive group \(\Q\). Let \(\mathcal D_\Gamma\) be the Hausdorff locally quasi-convex topologies on \(\Q\) with full character set \(\Gamma\) and continuous rational scalar maps. For \(\delta\in\mathcal D_\Gamma\), define
\[
T_f(\delta)=\init\{\mathrm{id}:\R\to\R_{\mathrm{usual}},\quad
f:\R\to(\Q,\delta)\}.
\]

\begin{theorem}[Transfer and incompatible joins]\label{q:thm:transfer-q}
The assignment \(T_f\) is an order embedding of \(\mathcal D_\Gamma\) into \(\mathcal C_{f,\Q}^{\Q}\), satisfies \(T_f(\gamma_0)=\sigma_{f,\Q}\), and preserves binary joins. Moreover
\[
T_f(\delta_1)\vee T_f(\delta_2)
\text{ is compatible with }H_{f,\Q}
\quad\Longleftrightarrow\quad
\widehat{(\Q,\delta_1\vee\delta_2)}=\Gamma.
\]
\end{theorem}
\begin{proof}
Rational linearity of \(f\) proves the TVGQ property. The initial construction gives Hausdorffness and local quasi-convexity. Put \(N=\ker f\). Every fiber of \(f\) is a translate of the usual-dense group \(N\). Consequently for a target neighborhood \(U\) and \(\varepsilon>0\),
\[
f\bigl((-\varepsilon,\varepsilon)\cap f^{-1}(U)\bigr)=U.
\]
Thus \(f:(\R,T_f(\delta))\to(\Q,\delta)\) is an open continuous surjection.

The induced topology on \(N\) is usual. Restrict any continuous character \(\chi\) to \(N\) and extend it to its usual completion \(\R\); the extension is \(p(cx)\). Then \(\chi-p(c\,\mathrm{id})\) vanishes on \(N\), and the open quotient map makes it a continuous character of \((\Q,\delta)\). Hence it is \(p(rf(x))\), \(r\in\Q\), giving exactly \(H_{f,\Q}\). The open quotient recovers \(\delta\) from \(T_f(\delta)\), proving order reflection. Order preservation, the formula for \(\gamma_0\), and preservation of joins follow from the generating maps.

If the rational join has dual \(\Gamma\), the preceding character calculation applies to it. Otherwise choose an extra continuous character \(\psi\). Its pullback \(\psi\circ f\) is continuous for the transferred join. If it belonged to \(H_{f,\Q}\), we would have \(\psi(f(x))=p(cx+rf(x))\). On the uncountable divisible group \(N\), this gives \(cx\in\Z\), hence \(cx=0\) by divisibility, and therefore \(c=0\). Surjectivity of \(f\) would then imply \(\psi(q)=p(rq)\), a contradiction. Thus the transferred join is incompatible.
\end{proof}

\begin{proposition}[A strict compatible refinement exists]\label{q:prop:strict-compatible}
There is \(\delta\in\mathcal D_\Gamma\) with \(\gamma_0\subsetneq\delta\). Consequently
\[
\sigma_{f,\Q}\subsetneq T_f(\delta)\in\mathcal C_{f,\Q}^{\Q}.
\]
In particular \(\sigma_{f,\Q}\) is not the greatest compatible TVGQ topology.
\end{proposition}
\begin{proof}
We use the established existence of a Hausdorff locally quasi-convex topology \(\lambda\) on \(\Q/\Z\), strictly finer than its circle subspace topology, with the same character group \(\Z\iota\), where \(\iota:\Q/\Z\to\T\) is the inclusion; see \cite{Barrera2014}. Give \(\Q\) the initial topology \(\delta\) of
\[
q\longmapsto rq+\Z\in(\Q/\Z,\lambda),\qquad r\in\Q.
\]
These maps separate points and are permuted by nonzero rational rescaling. The resulting topology is Hausdorff, locally quasi-convex, and rational-scalar-continuous.

A continuous character is controlled, by the quarter-arc criterion, by finitely many of these maps. Write their parameters as \(k_j/d\), with integers \(k_j\) and \(d>0\). Integer multiplication is \(\lambda\)-continuous, so these observations factor through the single surjection \(q\mapsto q/d+\Z\). The character is continuous for that single pullback topology. It vanishes on its kernel \(d\Z\), and therefore descends continuously to \((\Q/\Z,\lambda)\). Its parameter is thus an integer, giving on \(\Q\) a rational-parameter character. Conversely all characters in \(\Gamma\) are continuous, so \(\delta\in\mathcal D_\Gamma\).

The topology \(\lambda\) is not precompact: a precompact group topology is determined by its character set, so precompactness would make it equal to the original circle topology. The continuous surjection \(q\mapsto q+\Z\) shows that \(\delta\) is not precompact. Since \(\gamma_0\subseteq\delta\) and \(\gamma_0\) is precompact, the inclusion is strict. Theorem~\ref{q:thm:transfer-q} transfers strictness.
\end{proof}

Theorem~\ref{twochains:transfer} of Part I constructs an incompatible pair in $\mathcal D_\Gamma$ and disproves the former conjecture~\ref{q:conj:greatest-q}. In particular the cited failure of the usual topology on \(\Q\) to be Mackey does not itself supply that pair: its usual dual consists of all real-parameter characters, whereas the rational target here has the smaller dual \(\Gamma\). Nor is the transfer family proved cofinal among all compatible topologies on \(\R\).

We next derive a necessary quotient condition for greatestness and a Baire obstruction. These arguments apply to the rational Hamel pair \(H_{f,\Q}\), and do not assume that the conjectured topology exists unless explicitly stated.

\begin{lemma}[The rational quotient for every compatible topology]\label{q:lem:quotient-v4}
Let \(\eta\in\mathcal L(H_{f,\Q})\), without assuming rational scalar continuity, and put \(N=\ker f\). Then \(N\) is closed. Identify \(\R/N\) with \(\Q\) through \(f\), and let \(\lambda_\eta\) be the Hausdorff quotient group topology. Its full character group is exactly \(\Gamma\). The quotient map is an open continuous surjection, and \(\lambda_\eta\) is nondiscrete. Local quasi-convexity of the quotient is not asserted.
\end{lemma}
\begin{proof}
Every character \(p\circ qf\), \(q\in\Q\), belongs to the prescribed dual, and
\[
N=\bigcap_{q\in\Q}\ker(p\circ qf).
\]
Indeed, if \(f(x)\ne0\), choose \(q=1/(2f(x))\) to obtain the circle value \(1/2+\Z\). This proves closedness. A quotient of a topological group by a closed subgroup is Hausdorff, and its quotient map is open and continuous.

A character of this quotient pulls back to some \(p(cx+rf(x))\), \(c\in\R,r\in\Q\), vanishing on \(N\). Thus \(cN\subseteq\Z\). The group \(N\) is divisible, so an additive map from \(N\) into \(\Z\) is zero; hence \(cN=0\). Also \(N\ne\{0\}\), because \(f\) has rational rank one on the infinite-dimensional rational space \(\R\). Therefore \(c=0\). Conversely all \(p\circ rf\) vanish on \(N\) and descend continuously by the quotient definition, giving exactly \(\Gamma\).

If the quotient were discrete, every algebraic character of \(\Q\) would be continuous. For an irrational \(t\), the character \(q\mapsto p(tq)\) does not belong to \(\Gamma\): equality with a rational-parameter character would make the additive integer-valued map \(q\mapsto(t-r)q\) zero, forcing \(t=r\). This contradicts the dual calculation and proves nondiscreteness.
\end{proof}

\begin{theorem}[A Baire obstruction for the entire compatible class]\label{q:thm:nonbaire-v4}
Every topology in \(\mathcal L(H_{f,\Q})\), and in particular every compatible TVGQ topology in \(\mathcal C_{f,\Q}^{\Q}\), is non-Baire. Consequently no compatible member is locally compact or completely metrizable.
\end{theorem}
\begin{proof}
By Lemma~\ref{q:lem:quotient-v4}, its quotient by \(N\) is a countably infinite nondiscrete Hausdorff group. It has no isolated point, since translation of one isolated point would make the group discrete. Every singleton is therefore closed and nowhere dense. The countable union of these singletons is the whole group, so it is not Baire.

An open continuous surjection sends a Baire space to a Baire space. For completeness, inverse images of dense open sets are open and dense under such a map; applying the Baire property in each nonempty open inverse image proves density of the target intersection. Since the quotient map is open, a Baire source would give a Baire quotient, a contradiction. Locally compact Hausdorff spaces and completely metrizable spaces are Baire, proving the exclusions.
\end{proof}

\begin{theorem}[A necessary Mackey topology on the rational quotient]\label{q:thm:necessary-rational-v4}
If \(\mathcal C_{f,\Q}^{\Q}\) has a greatest member \(\mu\), then the dual pair \((\Q,\Gamma)\) has a greatest compatible Hausdorff locally quasi-convex group topology. This rational topology is also the greatest member of \(\mathcal D_\Gamma\).
\end{theorem}
\begin{proof}
Let \(\lambda=\lambda_\mu\) be the quotient topology of Lemma~\ref{q:lem:quotient-v4}. It has dual \(\Gamma\). Take any Hausdorff locally quasi-convex group topology \(\delta\) on \(\Q\) with dual \(\Gamma\). By Theorem~\ref{q:thm:saturation}, its rational saturation \(\delta^{[\Q]}\) belongs to \(\mathcal D_\Gamma\). Theorem~\ref{q:thm:transfer-q} and greatestness give
\[
T_f(\delta^{[\Q]})\subseteq\mu.
\]
Hence \(f:(\R,\mu)\to(\Q,\delta^{[\Q]})\) is continuous. The universal property of the quotient topology yields
\[
\delta\subseteq\delta^{[\Q]}\subseteq\lambda.
\]

Let \(\lambda_{\mathrm{qc}}\) be the join of all Hausdorff locally quasi-convex group topologies coarser than \(\lambda\). This family is nonempty: \(\gamma_0\) is one such topology because its defining characters belong to the dual of \(\lambda\). A join of locally quasi-convex group topologies is locally quasi-convex, by the diagonal product embedding. Therefore
\[
\gamma_0\subseteq\lambda_{\mathrm{qc}}\subseteq\lambda.
\]
Both endpoints have full dual \(\Gamma\). Monotonicity of continuous character groups gives the same dual for the middle topology. Every \(\delta\) considered above belongs to the family defining \(\lambda_{\mathrm{qc}}\), so is coarser than it. Thus \(\lambda_{\mathrm{qc}}\) is greatest among all compatible locally quasi-convex topologies for \((\Q,\Gamma)\). Corollary~\ref{q:cor:mackey-equivalence}, applied on \(\Q\), gives its rational scalar continuity and greatestness in \(\mathcal D_\Gamma\).
\end{proof}

\begin{corollary}[Restrictions on a proposed greatest topology]\label{q:cor:greatest-restrictions-v4}
If Conjecture~\ref{q:conj:greatest-q} holds, its greatest topology is non-Baire, nonprecompact, not locally compact, and not completely metrizable. The greatest rational topology forced by Theorem~\ref{q:thm:necessary-rational-v4} has the same four properties. In particular, finding a compatible locally compact or completely metrizable topology cannot prove the conjecture.
\end{corollary}
\begin{proof}
The Baire obstruction gives the source restrictions except nonprecompactness. A greatest source member contains \(\sigma_{f,\Q}\), which maps continuously onto the nonprecompact usual real line; it therefore cannot be precompact.

The rational topology is countable, Hausdorff, and nondiscrete by its dual \(\Gamma\), so the same singleton argument makes it non-Baire and excludes local compactness and complete metrizability. If it were precompact, its exact character set would force it to equal \(\gamma_0\). Proposition~\ref{q:prop:strict-compatible} supplies a strictly finer member of \(\mathcal D_\Gamma\), contradicting greatestness. Thus it is not precompact either.
\end{proof}

The necessary rational condition in Theorem~\ref{q:thm:necessary-rational-v4} fails: Theorem~\ref{twochains:transfer} of Part I proves that $(\Q,\Gamma)$ has no Mackey topology. More directly, its transferred pair disproves greatestness on the full rational Hamel pair. No assertion that every quotient of a locally quasi-convex group is locally quasi-convex is used.

\subsection{An obstruction for every dyadic partition}
The alternating-block attempt in \cite[Appendix A]{HamelCompanion} can be excluded for every partition, not only the particular blocks used there. For a nonempty set \(I\subseteq\mathbb N_0\), let \(\xi_I\) be the group topology on \(\Q\) generated by the seminorms
\[
Q_{I,d}(q)=\sup_{n\in I}\left\|\frac{q}{d2^n}\right\|_{\T},\qquad d\in\mathbb N_{>0}.
\]
Here bounded gaps means that some integers \(M,N\ge0\) make every interval \([k,k+M]\), \(k\ge N\), meet \(I\). These topologies are coarser than the usual topology because \(Q_{I,d}(q)\le |q|/d\).

\begin{proposition}[Bounded gaps and long blocks]\label{q:prop:dyadic-v4}
If \(I\) has bounded gaps, then \(\xi_I\) is the usual topology on \(\Q\). If \(I\) contains consecutive blocks of arbitrarily large length, then \(\xi_I\) has an irrational-parameter continuous character. Consequently, for every partition \(\mathbb N_0=I\cup J\), at least one of the two dyadic topologies has a character outside \(\Gamma\). No such partition gives two compatible members of \(\mathcal D_\Gamma\) with incompatible join.
\end{proposition}
\begin{proof}
Suppose first that \(I\) has bounded gaps, with constants \(M,N\). For \(|q|\ge2^{N-1}\), choose an integer \(k\ge N\) with \(1/4\le |q|/2^k<1/2\), and choose \(n\in I\cap[k,k+M]\). Then
\[
Q_{I,1}(q)\ge2^{-M-2}.
\]
Thus a sufficiently small \(Q_{I,1}\)-ball is contained in a fixed bounded interval \((-B,B)\). Fix \(n_0\in I\), and take an integer \(d\) with \(d2^{n_0}>2B\). On that interval, \(\|q/(d2^{n_0})\|_{\T}=|q|/(d2^{n_0})\). Intersecting the bounded ball with sufficiently small \(Q_{I,d}\)-balls therefore gives arbitrarily small usual neighborhoods. The reverse continuity follows from the estimate preceding the proposition.

Now suppose arbitrarily long blocks exist. Choose disjoint blocks \([a_j,b_j]\subseteq I\), \(j\ge1\), so that \(b_j-a_j\ge j\) and \(b_{j+1}-b_j\to\infty\). Such a choice is possible: arbitrarily long blocks remain arbitrarily far out, including when \(I\) contains a final infinite interval. Put
\[
\alpha=\sum_{j\ge1}2^{-b_j}.
\]
The endpoints can be chosen positive and strictly increasing. Its binary expansion has infinitely many ones and unbounded gaps, so \(\alpha\) is irrational. If \(Q_{I,1}(q)<\delta<1/4\), the centered representatives \(y_n\) of \(q/2^n+\Z\) satisfy \(y_n=2y_{n+1}\) within each block: the difference is an integer of absolute value less than \(3\delta<1\). Hence
\[
\left\|\frac{q}{2^{b_j}}\right\|_{\T}<\delta\,2^{-(b_j-a_j)},\qquad
\|\alpha q\|_{\T}\le\delta\sum_{j\ge1}2^{-(b_j-a_j)}\le\delta.
\]
The circle sums agree with the convergent real series modulo \(\Z\). This estimate for arbitrarily small \(\delta\) proves continuity of \(q\mapsto p(\alpha q)\), an extra character.

Finally, if \(I\) has bounded gaps its topology already has all real-parameter characters. Otherwise its complement \(J\) contains arbitrarily long consecutive blocks, and the second argument applies to \(J\). Empty parts are allowed here with the indiscrete topology; the full complementary part is usual. For any partition the join is usual, since the combined observations run over every nonnegative exponent, but individual compatibility already fails on at least one side.
\end{proof}

This rules out the dyadic partition route. It does not rule out observations from different chains: Theorems~\ref{twochains:H} and~\ref{twochains:transfer} of Part I construct the incompatible pair by that method. No uniform splitting lemma or converse to the rational-quotient reduction is needed.

\section{Conditions for a scalar field and a real vector topology}\label{q:sec:regularity-v3}
\begin{proposition}[An open mapping condition for inverse scalars]\label{q:prop:field-regularity}
If the additive group of a TVGQ is locally compact and sigma-compact, its continuous scalar ring is a field. Compact Hausdorff TVGQ spaces are included.
\end{proposition}
\begin{proof}
Every continuous surjective homomorphism \(\phi:G\to G\) of a locally compact sigma-compact group is open. To see this, for a zero-neighborhood \(U\) choose a relatively compact open neighborhood \(W\) with \(\overline W-\overline W\subseteq U\). Sigma-compactness supplies countably many translates of \(W\) covering \(G\). Their compact images \(\phi(\overline W)\), with the corresponding translates, cover the target. The Baire theorem for locally compact Hausdorff spaces implies that \(\phi(\overline W)\) has nonempty interior. Hence \(\phi(\overline W)-\phi(\overline W)\) contains a zero-neighborhood and is contained in \(\phi(U)\). Thus \(\phi\) is open. Apply this to each continuous bijection \(m_a\), \(a\ne0\). Its inverse \(m_{1/a}\) is continuous, so Proposition~\ref{q:prop:ring} makes the scalar ring a field.
\end{proof}

\begin{proposition}[Joint continuity from orbit continuity and uniform control]\label{q:prop:joint-v3}
Let \((V,\tau)\) be a TVGQ. Suppose each orbit map \(r\mapsto rx\), from the usual real line into \(V\), is continuous. Suppose also that the family \(\{m_q:q\in\Q,|q|\le1\}\) is equicontinuous at zero: every zero-neighborhood \(U\) contains \(qW\) for all these \(q\), for some zero-neighborhood \(W\). Then \(\tau\) is a real topological vector space topology, with jointly continuous scalar multiplication.
\end{proposition}
\begin{proof}
For a prescribed zero-neighborhood \(U\), choose a closed zero-neighborhood \(U_0\subseteq U\); Hausdorff topological groups are regular. Take \(W\) with \(qW\subseteq U_0\) for all rational \(|q|\le1\). For every real \(|r|\le1\) and \(w\in W\), approximate \(r\) by such rationals. Orbit continuity and closedness give \(rw\in U_0\). Thus all real scalars in the unit interval are equicontinuous at zero. Integer scaling extends this to every bounded real interval, proving continuity of each fixed real scalar map.

At \((r_0,x_0)\), use
\[
rx-r_0x_0=r(x-x_0)+(r-r_0)x_0.
\]
For \(r\) near \(r_0\), the first term is uniformly small when \(x\) is near \(x_0\), by bounded-interval equicontinuity. The second term is small by orbit continuity. Continuity of addition proves joint continuity.
\end{proof}

Even strong regularity of the group alone does not force every real scalar to act continuously. Choose a rational-linear bijection \(A:\R\to\R\) with \(A(1)=1\), \(A(\sqrt2)=2\sqrt2\), and fix all other elements of a Hamel basis containing \(1,\sqrt2\). Give \(\R\) the metric \(d(x,y)=|A(x)-A(y)|\). This TVGQ is additively isomorphic to the usual real line, and therefore connected, locally compact, sigma-compact, and completely metrizable. If \(m_{\sqrt2}\) were continuous, its conjugate \(A m_{\sqrt2}A^{-1}\) would be a usual-continuous additive map and hence multiplication by a constant. But it takes \(1\) to \(2\sqrt2\) and \(2\sqrt2\) to \(2\), contradicting that form. Its continuous scalar ring is a field by Proposition~\ref{q:prop:field-regularity}, but is a proper subfield of \(\R\).

\section{Isomorphism classification of the modular families}\label{q:sec:isomorphism-v3}
Use \(\mathcal F(f,S)=\{f_a:a\in S\}\) for the rational subspace of rational-valued additive functions supplied by the parameter space. The next theorem concerns arbitrary continuous additive homomorphisms, without assuming real linearity.

\begin{theorem}[Rigidity and exact isomorphism criterion]\label{q:thm:isomorphism-v3}
Let \(f,g\) be nonzero rational-valued rational-linear functionals and \(S,T\) nonzero rational subspaces. Fix one construction type \(\eta\), either \(\beta\) throughout or \(\sigma\) throughout. Every continuous additive homomorphism
\[
L:(\R,\eta_{f,S})\longrightarrow(\R,\eta_{g,T})
\]
has the form \(L(x)=cx\). For \(c\ne0\), this map is continuous exactly when \(\mathcal F(g,cT)\subseteq\mathcal F(f,S)\), and it is an isomorphism exactly when
\[
\mathcal F(g,cT)=\mathcal F(f,S).
\]
\end{theorem}
\begin{proof}
All ordinary characters \(p(rx)\) are target-continuous. With \(r=1\), Theorem~\ref{q:thm:characters} gives \(p(L(x))=p(dx+f_a(x))\), for some \(d\in\R,a\in S\). The additive integer-valued difference is zero, so \(L=d\,\mathrm{id}+f_a\). For a fixed irrational \(r\), the same argument gives \(rL=e\,\mathrm{id}+f_b\). Hence
\[
(rd-e)x=f_b(x)-rf_a(x).
\]
The right side has countable image, whereas any nonzero real multiple of the identity has uncountable image. Thus \(rd=e\) and \(rf_a=f_b\). Since both functional values are rational, irrationality forces \(f_a=0\), proving \(L=d\,\mathrm{id}\).

For \(L=m_c\), pullback of the target Hamel character indexed by \(a\in T\) is \(p\circ g_{ca}\). Its membership in \(H_{f,S}\) is equivalent to \(g_{ca}=f_b\) for some \(b\in S\): the possible ordinary term vanishes by rational-valuedness and Lemma~\ref{q:lem:mixed}'s integer-valued argument. This gives the necessary inclusion. It also suffices by the initial-topology definitions; in the \(\sigma\) case the usual coordinate is already continuous. Applying the same reasoning to \(m_{1/c}\) gives the equality criterion for an isomorphism.
\end{proof}

\begin{corollary}[Consequences after forgetting the real action]\label{q:cor:isomorphism-v3}
For fixed \(f\), the spaces \((\R,\eta_{f,S})\) and \((\R,\eta_{f,T})\), for either modular construction type, are additively topologically isomorphic exactly when \(S=cT\) for some \(c\ne0\). Within each modular family the continuous scalar ring is an invariant even of additive topological group isomorphism. In particular distinct intermediate fields used as parameter spaces give nonisomorphic groups, regardless of the choice of functional. For \(S=T=\Q\), the isomorphism criterion is
\[
g(cx)=r f(x)\quad(x\in\R),\qquad c\in\R\setminus\{0\},\quad r\in\Q\setminus\{0\}.
\]
Each of the families \(\{\beta_{f,\Q}:f\ne0\}\) and \(\{\sigma_{f,\Q}:f\ne0\}\) contains exactly \(2^{\mathfrak c}\) additive topological group isomorphism types.
\end{corollary}
\begin{proof}
Injectivity of \(a\mapsto f_a\) gives the first assertion. Every isomorphism furnished by Theorem~\ref{q:thm:isomorphism-v3} is a real scalar map, so commutes with every scalar multiplication. Conjugation therefore preserves exactly which scalar maps are continuous. For a parameter field the scalar ring is that field, proving the distinction. When \(S=T=\Q\), equality of the one-dimensional rational observation spaces is precisely the displayed condition.

There are \(2^{\mathfrak c}\) nonzero rational-linear functionals: arbitrary rational assignments on a Hamel basis of size \(\mathfrak c\) give exactly that cardinality. For fixed \(f\), the displayed isomorphism condition permits at most \(\mathfrak c\) choices of \(g\), because \(g(x)=r f(x/c)\). Partitioning \(2^{\mathfrak c}\) elements into classes of size at most \(\mathfrak c\) leaves exactly \(2^{\mathfrak c}\) classes. This proves the count for either construction type.
\end{proof}

The two modular types cannot be isomorphic to one another: $\beta_{f,S}$ is precompact, whereas $\sigma_{g,T}$ has the usual real line as a continuous quotient and is not precompact. Section~\ref{q:sec:tauiso} gives the separate classification for real-valued observation spaces, where nonscalar automorphisms occur. Scalar rings are not invariants of arbitrary topological groups with unspecified real actions, as the transported real-space example in Section~\ref{q:sec:regularity-v3} shows.

\section{Scope and further questions}
Theorem~\ref{q:thm:universal-complete} settles existence of a universal separating sequence for every parameter space: it is equivalent to a finite scalar covering of \(\R\) by \(S\). This includes the formerly unaddressed proper subspaces of dimension \(\mathfrak c\). The rational case is negative, every fixed functional still has a separating sequence, and the unscaled Pell sequence has the sharper particular criterion of Theorem~\ref{q:thm:pell-criterion}.

For the greatest-compatible-topology question, compatible rational saturation identifies greatestness with the unrestricted Mackey problem for the same dual pair. Theorems~\ref{twochains:H} and~\ref{twochains:transfer} of Part I give two compatible topologies with incompatible join and disprove greatestness for every nonzero $f$. The same quotient pair also transfers to the full real-orbit pair, disproving the real-scalar assertion. All compatible topologies for the rational Hamel pair remain non-Baire by Theorem~\ref{q:thm:nonbaire-v4}.

Local compactness and sigma-compactness force the continuous scalar ring to be a field, and orbit continuity together with rational equicontinuity forces joint real scalar continuity. Group regularity alone does not imply the latter. Theorem~\ref{q:thm:isomorphism-v3} classifies additive isomorphisms within each modular family and proves that its scalar-ring distinctions survive forgetting the real action. Section~\ref{q:sec:tauiso} recalls the cited classification of arbitrary additive isomorphisms for the real-valued spaces $\tau_{f,S}$.

The remaining problems are unresolved in this article; this is not an assertion that the literature leaves them open. The operator-group framework and compact duality are established tools; the present article gives direct proofs for the specified Hamel families.

Section~\ref{q:sec:chains-audited} recalls the cited character criteria for rational and torsion observation topologies from divisibility chains. Compatible subsets of one chain have compatible union, so that family cannot yield the required incompatible pair. Proposition~\ref{chains:irrational} of Part I gives sequence witnesses for each excluded irrational scalar, and Proposition~\ref{chains:probe} of Part I gives a fixed-character join test for the rational Hamel pair. The counterexample of Section~\ref{sec:twochains} of Part I settles both greatest-compatible existence and the general compact-union question negatively; Section~\ref{q:sec:twochains} explains its TVGQ application.

\subsection*{Source comparison and an exclusion for coarse topologies}
Kuba \cite{Kuba2026} studies Hausdorff topologies on the real line coarser than the usual topology, including counting results up to homeomorphism. These results do not assert continuity of addition or of fixed rational scalar multiplications, or preservation of the character group $H_{f,\Q}$. The following elementary argument excludes one of his classes as a source of compatible TVGQ counterexamples.

\textbf{Exclusion for the class $L_0$.}
Let $\tau$ be a Hausdorff group topology on the additive real line with $\tau\subseteq\tu$. If its restriction to $\R\setminus\{0\}$ is the usual subspace topology, then $\tau=\tu$.

\begin{proof}
Hausdorffness makes $\{0\}$ closed, so its complement is $\tau$-open. The restriction hypothesis therefore makes every usual-open interval avoiding zero $\tau$-open. For every $\varepsilon>0$, translate $(1,1+\varepsilon)$ by $-1-\frac{\varepsilon}{2}$. Translation is a homeomorphism for an additive group topology, so $(-\frac{\varepsilon}{2},\frac{\varepsilon}{2})$ is $\tau$-open. Further translations show that every usual-open set is $\tau$-open. The assumed reverse inclusion proves equality.
\end{proof}

Kuba's class $L_0$ consists of Hausdorff coarsenings with this restriction property. Its strict coarsenings are consequently not additive group topologies and hence are not TVGQ topologies. The proof does not require continuity of any scalar multiplication, so the exclusion applies equally to TVGQ and real TVG topologies. Moreover the usual real-line topology itself is not compatible with $H_{f,\Q}$: it does not make the nonzero Hamel character $p\circ f$ continuous, by Lemma~\ref{q:lem:usual}. Thus $L_0$ contributes no member to the compatible rational Hamel class. This statement concerns $L_0$, not every topology in \cite{Kuba2026} or every possible adaptation of its methods. No historical novelty is claimed for the exclusion argument, which is also recorded in the companion article \cite{HamelCompanion}.

Hejcman's paper on boundedness \cite{H1959} is a separate work from his 1988 paper on real-line TVG topologies \cite{H1988}. The available archive listings and conference report do not establish the theorem content of the latter; no comparison with its full text is claimed here. These sources do not establish or refute Conjecture~\ref{q:conj:greatest-q}, do not supply the incompatible pair required by Theorem~\ref{q:thm:greatest-criterion}, and do not verify historical priority for the specialized results of this article. Theorem~\ref{twochains:transfer} of Part I disproves the conjecture by a separate argument. The real-valued isomorphism classification is recalled in Section~\ref{q:sec:tauiso} from Section~\ref{sec:tauiso} of Part I.

\section{Divisibility-chain results used here}\label{q:sec:chains-audited}

The rational and torsion chain criteria, with their proofs, are established in Section~\ref{sec:chains-audited} of Part I. We recall the scope of those results to distinguish the same-chain obstruction from the two-chain counterexample.

For an increasing divisibility chain, compatibility holds exactly when, for each fixed bound, the lengths of consecutive runs of ratios below that bound are uniformly bounded. On $\Q$, this is the criterion for the full dual to be $\Gamma$ for the rational observation topologies with all positive denominator parameters. On $\Q/\Z$, it is the criterion for the full dual to be $\Z\iota$, where $\iota$ is the natural inclusion into $\T$; the infinite torsion chain starts with $1$ and every positive integer must divide some term. See Theorems~\ref{chains:F} and~\ref{chains:G} of Part I.

Compatible subsets of one chain have a compatible union (Corollaries~\ref{chains:Fjoin} and~\ref{chains:Gjoin} of Part I). Thus splitting a single chain cannot produce the incompatible pair needed here. The factorial-chain null sequence, excluded irrational-scalar witnesses and fixed-character join test are given in Propositions~\ref{chains:factorial}, \ref{chains:irrational} and~\ref{chains:probe} of Part I. These results are used as cited input, rather than proved a second time.

\section{The two-chain counterexample in the TVGQ framework}\label{q:sec:twochains}

Theorem~\ref{twochains:H} of Part I constructs two compatible topologies $\lambda_1,\lambda_2$ on $\Q/\Z$, with full dual $\Z\iota$, whose join has an extra character. The construction and its Chinese remainder and character arguments belong to Part I. Theorem~\ref{twochains:transfer} there supplies their rational lifts $\delta_1,\delta_2\in\mathcal D_\Gamma$ with an incompatible join.

For clarity, with $\pi(q)=q+\Z$, these lifts are

\[
\delta_i=\operatorname{init}\{q\mapsto\pi(rq):\Q\to(\Q/\Z,\lambda_i),\ r\in\Q\},\qquad i=1,2.
\]

Their full dual is $\Gamma$, every rational dilation is continuous, and some character $\psi\notin\Gamma$ is continuous on $\delta_1\vee\delta_2$. These are the rational conclusions established in the proof of Theorem~\ref{twochains:transfer} of Part I.

We now apply the rational transfer already proved in Theorem~\ref{q:thm:transfer-q}. For each nonzero rational-linear $f:\R\to\Q$, put

\[
\eta_i=T_f(\delta_i)=\operatorname{init}\{\mathrm{id}:\R\to\R_{\mathrm{usual}},\ f:\R\to(\Q,\delta_i)\}.
\]

Then $\eta_1,\eta_2$ are compatible Hausdorff locally quasi-convex TVGQ topologies with full dual $H_{f,\Q}$. The join criterion in Theorem~\ref{q:thm:transfer-q} makes $\eta_1\vee\eta_2$ incompatible. Theorem~\ref{q:thm:greatest-criterion} therefore rules out a greatest compatible TVGQ topology and an unrestricted Mackey topology for $(\R,H_{f,\Q})$. This is the application of the primary counterexample to the rational-scalar framework, not a second counterexample.

The full real-scalar transfer and the nonexistence conclusions for the rational and torsion dual pairs are proved in Theorem~\ref{twochains:transfer} of Part I. The common quotient counterexample answers both scalar questions, but the all-real-coordinate transfer is not needed for the TVGQ application above.

\section{Real-valued isomorphism results used here}\label{q:sec:tauiso}

The full homomorphism and isomorphism classification for the real-valued observation spaces is established in Section~\ref{sec:tauiso} of Part I. It already allows arbitrary rational parameter spaces, so no additional rational-case proof is required here. Write $\mathcal F(f,S)=\{x\mapsto f(ax):a\in S\}$.

For nonzero rational-linear $f,g:\R\to\Q$ and rational subspaces $S,T$, an additive map $L:(\R,\tau_{f,S})\to(\R,\tau_{g,T})$ is continuous exactly when it has the form

\[
L(x)=cx+\sum_{j=1}^{n}d_j f(b_jx),\qquad c,d_j\in\R,\quad b_j\in S,
\]

with $\mathcal F(g,cT)\subseteq\mathcal F(f,S)$ when $c\ne0$; when $c=0$ no further condition is needed (Theorem~\ref{tauiso:hom} of Part I). The invertibility criterion and explicit inverse for such maps are given in Theorem~\ref{tauiso:iso} of Part I.

In particular, an additive topological isomorphism exists if and only if a scalar isomorphism exists, equivalently

\[
\mathcal F(g,cT)=\mathcal F(f,S)\qquad\text{for some }c\in\R\setminus\{0\}.
\]

The continuous scalar ring is consequently an isomorphism invariant within this family, and the rational-parameter family has exactly $2^{\mathfrak c}$ isomorphism types (Corollary~\ref{tauiso:types} of Part I). Individual isomorphisms need not be scalar.

For $S=\Q$, the automorphisms are exactly

\[
L_{c,s}(x)=cx+s f(x),\qquad c\in\Q\setminus\{0\},\quad s\in\R,\quad c+f(s)\ne0,
\]

as proved in Corollary~\ref{tauiso:auto} of Part I. That result also gives the composition and inverse formulas and the nonabelian shear subgroup. We use these conclusions in the comparison below without repeating their proofs.

\section{Comparison with the full-parameter (TVG) case}\label{q:sec:comparison}

Every TVG is a TVGQ, whereas the converse fails. In the Hamel constructions, the full real parameter space $S=\R$ gives the TVG case; a nonzero rational subspace $S\subseteq\R$ gives the general TVGQ family. The comparison below concerns these constructions and their specified compatible classes, rather than arbitrary topological groups.

Fix a nonzero rational-linear $f:\R\to\Q$. Write $M(S)=\{r\in\R:rS\subseteq S\}$, $K_S=\widehat{S_d}$, and $H_{f,S}=\{x\mapsto p(cx+f(ax)):c\in\R,\ a\in S\}$. Here $S_d$ is the discrete additive group, the hat denotes its Pontryagin dual, and $B$ denotes a rational basis of $S$. The topology $\beta_{f,S}$ is generated by $H_{f,S}$, $\sigma_{f,S}=\tu\vee\beta_{f,S}$, and $\tau_{f,S}$ is defined in Section~\ref{q:sec:three}. Subscripts $S=\R$ are omitted.

A shared assertion below means the same structural statement with the appropriate parameter range. It does not identify the prescribed dual groups. The TVG column is the case $S=\R$ of the results of this part, except for the rows on precompact classification, the two-chain counterexample and the greatest compatible topology, and the real-valued homomorphisms and isomorphisms, for which the full real-orbit statements are established in \cite{HamelCompanion} (Sections 9, 10, 13 and 15 for the modular, completion, rigidity and classification statements; Section 14 for the sequence statements) and in Part I (Section~\ref{sec:twochains} for the two-chain counterexample and its real transfer, Section~\ref{sec:tauiso} for the real-valued homomorphisms and isomorphisms, Section~\ref{sec:sequences} for the universal-sequence statements).

\begingroup
\small\setlength{\tabcolsep}{4pt}\renewcommand{\arraystretch}{1.15}
\begin{longtable}{>{\raggedright\arraybackslash}p{0.22\textwidth}>{\raggedright\arraybackslash}p{0.34\textwidth}>{\raggedright\arraybackslash}p{0.34\textwidth}}
\caption{Structural results for the Hamel families}\\
\hline
\textbf{Statement} & \textbf{TVG: $S=\R$} & \textbf{TVGQ: general $S$} \\ \hline
\endfirsthead
\hline
\textbf{Statement} & \textbf{TVG: $S=\R$} & \textbf{TVGQ: general $S$} \\ \hline
\endhead
\hline\endfoot
Scalar continuity (different) & Every fixed real multiplication is continuous. & Every fixed rational multiplication is continuous. Neither definition requires joint scalar continuity. \\ \hline
Exact scalar ring (different) & $C_\eta=\R$ for $\eta=\beta_f,\sigma_f,\tau_f$. & $C_\eta=M(S)$ for all three constructions. For $S=\Q$, this is exactly $\Q$. \\ \hline
Inverse scalar maps (different) & Every nonzero scalar multiplication is a homeomorphism. & For $r\ne0$, it is a homeomorphism exactly when $rS=S$. A continuous scalar bijection can have discontinuous inverse. \\ \hline
Basic properties (same) & All three topologies are Hausdorff and nondiscrete; $\beta_f$ is precompact. & The same holds for every nonzero $S$. \\ \hline
Strict separation (same) & $\tu\subsetneq\sigma_f\subsetneq\tau_f$. & $\tu\subsetneq\sigma_{f,S}\subsetneq\tau_{f,S}$. \\ \hline
Modular characters (same form) & Both modular topologies have full dual $H_f=H_{f,\R}$. & Both have full dual $H_{f,S}$. The observation parameter is restricted to $S$. \\ \hline
Real-valued additive maps (same) & On $\beta_f$, only zero; on $\sigma_f$, exactly $x\mapsto cx$. & The same conclusions for $\beta_{f,S}$ and $\sigma_{f,S}$. \\ \hline
Modular completions (same form) & For $\beta_f$: $K_{\R}\times K_{\R}$. For $\sigma_f$: $\R\times K_{\R}$. & For $\beta_{f,S}$: $K_{\R}\times K_S$. For $\sigma_{f,S}$: $\R\times K_S$. \\ \hline
Real-valued completion (same form) & $\R\times\R^B$, with $B$ a rational basis of $\R$. & $\R\times\R^B$, with $B$ a rational basis of $S$. For $S=\Q$, this is $\R^2$. \\ \hline
Modular rigidity (same) & Continuous additive maps within either modular family are scalar multiplications. & The same rigidity holds. For $c\ne0$, continuity is equivalent to $\mathcal F(g,cT)\subseteq\mathcal F(f,S)$. \\ \hline
Real-valued homomorphisms (same form) & Finite sums $cx+\sum_j d_jf(b_jx)$, with the target-observation condition. & The same description, with $b_j\in S$. Section~\ref{q:sec:tauiso} gives the precise condition. Nonscalar shears occur when $1\in S$. \\ \hline
Real-valued isomorphisms (same) & An isomorphism exists exactly when a scalar isomorphism exists. & The same: $\mathcal F(g,cT)=\mathcal F(f,S)$ for some $c\ne0$. Individual isomorphisms can be nonscalar. \\ \hline
Precompact classification (different) & Nonzero real vector subspaces of the algebraic character space; separation is automatic. & Point-separating rational vector subspaces; separation must be required. Scalars act by precomposition. \\ \hline
Two-chain counterexample (same source) & Transfers to compatible TVGs for $H_f$ with incompatible join. & Transfers to compatible TVGQs for $H_{f,\Q}$ with incompatible join. The rational source is the same; the transfers differ. \\ \hline
Greatest compatible topology (same negative answer) & No greatest compatible TVG for $H_f$, and no unrestricted Mackey topology for this dual pair. & The corresponding conclusions hold for $H_{f,\Q}$. This counterexample alone does not settle every general $S$. \\ \hline
\end{longtable}
\endgroup

\begingroup
\small\setlength{\tabcolsep}{4pt}\renewcommand{\arraystretch}{1.15}
\begin{longtable}{>{\raggedright\arraybackslash}p{0.22\textwidth}>{\raggedright\arraybackslash}p{0.34\textwidth}>{\raggedright\arraybackslash}p{0.34\textwidth}}
\caption{Sequence separation and its quantifiers}\\
\hline
\textbf{Statement} & \textbf{TVG: $S=\R$} & \textbf{TVGQ: general $S$} \\ \hline
\endfirsthead
\hline
\textbf{Statement} & \textbf{TVG: $S=\R$} & \textbf{TVGQ: general $S$} \\ \hline
\endhead
\hline\endfoot
Separation for each fixed $f$ & A separating sequence exists. & A separating sequence exists for every nonzero $S$; it may depend on $f$. \\ \hline
One sequence for every $f$ & The unscaled Pell sequence works. & Such a sequence exists exactly when $b_1S+\cdots+b_mS=\R$ for finitely many real $b_j$. \\ \hline
Unscaled Pell sequence & Works for every nonzero $f$. & Works for every nonzero $f$ exactly when $S+\sqrt2S=\R$. \\ \hline
Rational parameter case & $S=\Q$ is not the full TVG case. & No universal separating sequence for $S=\Q$, although each fixed $f$ has its own. \\ \hline
Universal modular convergence & Convergence for every $f$ forces finite rational span. & The same necessity holds for every nonzero $S$. \\ \hline
\end{longtable}
\endgroup

For the sequence table, separation means convergence to zero in $\sigma_{f,S}$ (and also in $\beta_{f,S}$ for the sequences constructed) but failure of convergence to zero in $\tau_{f,S}$. The unscaled Pell sequence is $z_n=n!(3-2\sqrt2)^{n^2}$. It is modularly null for every $f,S$; failure of real-valued convergence occurs exactly when $f$ is nonzero on $S+\sqrt2S$. Thus a universal separator and a separator chosen after fixing $f$ are different assertions.

The finite-covering criterion explains the distinction: universal modular convergence forces the sequence into a finite rational span $W$. If the rational span of $SW$ is proper, a nonzero functional vanishing on it makes every observation zero. Conversely, a finite covering by scalar multiples of $S$ permits finitely many factorially scaled rational approximations to irrational numbers to be interleaved; their observations are eventually integral, while at least one real-valued observation fails to tend to zero for each nonzero $f$. For $S=\Q$, no finite covering of $\R$ is possible.

These tables compare mathematical statements rather than levels of formal verification; the verification status is recorded in Section~\ref{q:sec:further-results}.

\section{Further results and formal verification}\label{q:sec:further-results}

Part I gives further results for the rational and torsion constructions used here, and for the full parameter space $S=\R$. We record their consequences and limits. Part I is the primary source for the chain criteria, the two-chain counterexample and the real-valued isomorphism classification. Sections~\ref{q:sec:chains-audited}--\ref{q:sec:tauiso} above recall these results and explain their use in the rational-scalar framework; their proofs are confined to Part I. Section~\ref{q:sec:comparison} compares the full-parameter and general-parameter cases.

\subsection{Nuclearity of the rational and quotient-chain topologies}

Write $\ell\in\{2,3\}$ for the prime and retain $p:\R\to\T$ for the quotient map. Define

\[
P_{\ell,d}(q)=\sup_{n\ge0}\left\|p\left(\frac{q}{d\ell^{n^2}}\right)\right\|_{\T},\qquad q\in\Q,\quad d\ge1.
\]

Let $\gamma_\ell$ be generated by the balls $P_{\ell,d}<\varepsilon$. Theorem~\ref{chains:F} of Part I gives its full dual $\Gamma$, because the consecutive ratios of $\ell^{n^2}$ tend to infinity. Both $(\Q,\gamma_2)$ and $(\Q,\gamma_3)$ are non-nuclear, and both quotient-chain groups $(\Q/\Z,\lambda_{F_1})$ and $(\Q/\Z,\lambda_{F_2})$ of Theorem~\ref{twochains:H} of Part I are non-nuclear; see Section~\ref{sec:nuclearity} of Part I. Thus nuclearity is no longer an unresolved issue for these four groups.

The obstruction uses the established theorem that nuclear groups respect compactness. If a Hausdorff group has a countable point-separating full character group $E$ and respects compactness, its topology equals the metrizable weak topology $\sigma(G,E)$. Indeed every weakly convergent sequence, together with its limit, is compact in both topologies; the compact-to-Hausdorff identity on that set is a homeomorphism. Sequential continuity of the identity from the metrizable weak topology proves equality. In particular such a nuclear group must be precompact.

For the prime groups a concrete witness is

\[
a_n^{(\ell)}=\frac{n!}{\ell^{v_\ell(n!)}}\,\ell^{n^2-1},\qquad n\ge1.
\]

Here $v_\ell(n!)$ is the exponent of $\ell$ in $n!$. Every fixed positive integer divides $a_n^{(\ell)}$ eventually, so this sequence tends to zero in $\sigma(\Q,\Gamma)$. However $P_{\ell,1}(a_n^{(\ell)})\ge\frac{1}{\ell}$ for $\ell=2,3$. For the quotient-chain groups $(\Q/\Z,\lambda_{F_i})$ the full dual is $\Z\iota$, which is countable and point-separating, so a precompact topology with this dual would be the circle subspace topology. In the circle topology no infinite family of distinct characters $n\iota$, $n\in N\subseteq\Z$, is equicontinuous. Indeed, given $\delta>0$, choose $n\in N\setminus\{0\}$ with $\frac{1}{2|n|}<\delta$. Then $x=\frac{1}{2|n|}+\Z$ lies in the $\delta$-neighbourhood of zero while $\|nx\|_{\T}=\frac{1}{2}$. The defining observations $x\mapsto M_jx+\Z$, with distinct $M_j\to\infty$, are equicontinuous for $\lambda_{F_i}$ because the uniform-convergence balls
\[
U_\varepsilon=\left\{x\in\Q/\Z:\sup_j\|M_jx\|_{\T}<\varepsilon\right\},\qquad \varepsilon>0,
\]
are zero-neighbourhoods. Hence $\lambda_{F_i}$ is not precompact, and the nuclearity argument above excludes nuclearity. Part I (Section~\ref{sec:nuclearity}) supplies the nuclear-group reference and the full arguments. These conclusions concern the rational and torsion groups themselves; they do not by themselves establish non-nuclearity of every real-line transfer.

\subsection{Sequence results and the parameter space}

For the full parameter space $S=\R$, Section~\ref{sec:sequences} of Part I proves an exact universal convergence criterion. Write $\tau_f=\tau_{f,\R}$ and similarly for $\sigma_f,\beta_f$. A sequence converges to zero in $\tau_f$ for every nonzero rational-linear $f:\R\to\Q$ if and only if its rational span is finite-dimensional and its rational coordinates tend to zero usually.

For convergence to zero in every $\sigma_f$, equivalently every $\beta_f$, the criterion is finite-dimensional rational span together with

\[
x_n=\sum_{j=1}^{m}q_{n,j}b_j\longrightarrow0\ \text{usually},\qquad \left\|p\left(\frac{q_{n,j}}d\right)\right\|_{\T}\longrightarrow0\quad(j=1,\ldots,m,\ d\ge1),
\]

where $b_1,\ldots,b_m$ is a rational basis of the span. In finite rational dimension these tests also apply to each fixed nonzero $f$ in the full-parameter family. Infinite-dimensional behavior can depend on the functional: there exist nonzero $f,g$ and a rationally independent sequence converging to zero in $\tau_f$ but not in any of $\tau_g,\sigma_g,\beta_g$.

These full-parameter statements should not be substituted for an individual-functional coordinate criterion with an arbitrary smaller $S$. The universal separating-sequence theorem of Section~\ref{q:sec:sequences} already handles every nonzero parameter space $S$, with its distinct quantifiers. The new results do not classify all null sequences for an arbitrary fixed functional with infinite-dimensional rational span.

\subsection{Size of the mixed-prime dual}

Put $P(k)=\max\{P_{2,1}(k),P_{3,1}(k)\}$ on $\Z$ and let

\[
H=\{\alpha+\Z\in\T:k\mapsto p(\alpha k)\text{ is continuous on }(\Z,P)\}.
\]

Theorem~\ref{mixed:dimension} of Part I and the remarks following it prove that $H$ has Hausdorff dimension zero in the usual circle metric and has Lebesgue measure zero, while $\Z[1/6]/\Z\subseteq H$. More quantitatively, for $0<\varepsilon<1/4$ the compact polar

\[
K_\varepsilon=\{\alpha+\Z:\|p(\alpha k)\|_{\T}\le1/4\text{ whenever }P(k)<\varepsilon\}
\]

can be covered by at most $\exp(C_\varepsilon\sqrt{\log N})$ circle balls of radius $1/(50N)$, for every integer $N\ge2$, with $C_\varepsilon$ independent of $N$. Consequently each such polar has upper box dimension zero.

The dimension estimate alone does not establish whether $H\subseteq\Q/\Z$. A dimension-zero set can contain irrational parameters, so neither this estimate nor non-nuclearity calculates the full character group of $\gamma_2\vee\gamma_3$. This is the specific mixed-prime problem of \cite[Section 22]{HamelCompanion}, distinct from the two-chain counterexample.

\subsection{Checks by Lean and Mathlib}

Selected statements were formalized in Lean 4 with Mathlib and typechecked. The checked declarations depend only on the standard axioms \texttt{propext}, \texttt{Classical.choice} and \texttt{Quot.sound}. A formal statement is checked under its stated hypotheses (for the Hamel functional it is normalized by a point $u$ with $f(u)=1$, which is no restriction). Formal verification does not certify every printed proof or the wording of the article.

Checked: the division and homeomorphism criteria for rational scalars; the finite-coordinate lemmas; the completions of $\sigma_{f,S}$, $\beta_{f,S}$ and $\tau_{f,S}$; the full modular character groups and the real-valued additive maps; strict separation; the exact continuous scalar rings and the order theorem; the rational case; the classification of precompact TVGQ topologies and the minimal rational topology; the Pell and scaled sequence results, the obstruction to a universal sequence, the complete criterion for a universal sequence and its field corollary via Artin--Schreier; directed joins, rational scalar saturation and its equivalence with the Mackey question, the scalar restrictions and the greatest-member criterion; the transfer theorem with its dual and join criterion; the strict compatible refinement (formalized with an explicit factorial-chain topology in place of the cited existence theorem); the rational quotient, the Baire obstruction and the necessary Mackey condition on the rational quotient; the restrictions on a greatest topology; the dyadic partition obstruction; the field and joint-continuity propositions; the rigidity and isomorphism results for the modular families and the counting of isomorphism types; the chain criteria, the two-chain construction and its rational and real transfers, and the homomorphism, isomorphism, type and automorphism results for the real-valued observation spaces proved in Part I; the non-precompactness and metrizability of $\gamma_\ell$, the explicit weakly null sequence and its failure of convergence in the prime topology, and, for the mixed-prime character parameters, the quantitative covering estimate, Hausdorff dimension zero, Lebesgue measure zero and the rational-subgroup inclusion. The formal proof of the covering estimate uses discrete grid orthogonality in place of the integral Parseval argument.

Not formalized: the theory of nuclear groups and the non-nuclearity conclusions that rely on it; the cited existence theorem of \cite{Barrera2014} (replaced by the explicit chain topology mentioned above); the upper box dimension conclusion and the reverse inclusion in the polar-union description; the exclusion argument for Kuba's class $L_0$; the comparison tables; and all bibliographic statements. No historical priority claim is made.

The Lean and Mathlib sources for the selected formalized results are available in the public repository \href{https://github.com/irinayar/MyPseudotopology-Lean}{MyPseudotopology-Lean}. The source snapshot cited here is \href{https://github.com/irinayar/MyPseudotopology-Lean/tree/c0444e290d5b7f0952907356f7f6a23e62b87586}{commit \texttt{c0444e2}}. Its verification records state the checked declarations and their source hashes. Individual module checks are distinguished from a full rebuild of the entire library.

\section*{Acknowledgement}
The author acknowledges OpenAI's Codex and Anthropic's Claude for assistance with mathematical discussion, review of arguments, exposition, and article preparation. The author also gratefully acknowledges the developers and contributors of the Lean proof assistant and the mathlib library, which were used for formal verification of selected results. The author assumes responsibility for the final content and any remaining errors.

\end{document}